\documentclass[11pt,draftcls,onecolumn]{IEEEtran}
\usepackage{etex} 
\usepackage{amsmath,graphicx,amssymb,epsfig,setspace,multicol,cite}
\usepackage{multirow}
\usepackage{float,dsfont}
\usepackage[latin1]{inputenc}
\usepackage{graphicx}
\usepackage{mathtools}
\usepackage{ulem} 
\usepackage{bbm} 
\usepackage{schemabloc}  
\usepackage{epsfig}       
\usepackage{pst-all}      
\usepackage{cancel}

\usepackage{tikz}
\usepackage{epstopdf}
\usetikzlibrary{positioning}

\DeclareTextSymbol{\degre}{T1}{6}
\DeclareTextSymbol{\degre}{OT1}{23}

\newcommand\tra{{\rm Tr}}

\newcommand\cov{{\rm Cov}}
\newcommand\var{{\rm Var}}
\newcommand\ve{{\rm vec}}
\newcommand\diag{{\rm Diag}}
\newcommand\e{{\rm E}}
\newcommand\pardef{ \stackrel{{\rm def}}{=} }

\newtheorem{definition}{Definition}
\newtheorem{theorem}{Theorem}

\usepackage{soul} 
\usepackage{color}

\begin{document}
\title{Background on real and complex elliptically symmetric distributions
\\
\large{Chapter published in the book "Elliptically symmetric distributions in Signal Processing and Machine Learning"by Springer with corrections and complements}}

\author{Jean-Pierre Delmas
\thanks{
Jean-Pierre Delmas is with Samovar laboratory, Telecom SudParis, Institut Polytechnique de Paris, 91120 Palaiseau, France,
e-mail: jean-pierre.delmas@it-sudparis.eu.
}}

\maketitle

\begin{abstract}
This chapter presents a short overview of real elliptically symmetric (RES) distributions, complemented by  circular complex elliptically symmetric (C-CES) and noncircular CES (NC-CES) distributions as complex representations 
of RES distributions. 
These distributions are both an extension of the multivariate Gaussian distribution and a  multivariate extension of univariate symmetric distributions.
They are equivalently defined through their characteristic functions and their stochastic representations, which naturally follow from the spherically symmetric distributions after affine transformations.
Particular attention is paid to the absolutely continuous case and to the subclass of compound Gaussian distributions. 
Results related to moments, affine transformations, marginal and conditional distributions, and summation stability are also presented.
Some well-known instances of RES distributions are provided with their main properties. 
Finally, the estimation of the symmetry center and scatter matrix is briefly discussed through the sample mean (SM), sample covariance matrix (SCM) estimate, maximum estimate (ML), $M$-estimators, and Tyler's $M$-estimators. Particular attention will be paid to the asymptotic Gaussianity of the $M$-estimators of the scatter matrix. 
To conclude, some hints about the Slepian-Bangs formula are provided.
\end{abstract}
\section{Introduction}
\label{sec:Introduction}
%
Until fifty years ago, most of the procedures in multivariate analysis were developed under the Gaussian assumption, mainly for mathematical convenience. However, in many applications, Gaussianity is a poor approximation of reality.  As a consequence, elliptically symmetric distributions have been widely used in various applications due to their 
flexibility and capability to better model various data behavior.
These distributions form a natural extension of the Gaussian one by allowing for both heavier-than-Gaussian and lighter-than-Gaussian tails while maintaining the elliptical geometry of the underlying equidensity (when it exists) contours.
These real elliptically symmetric (RES) distributions were equivalently defined in the statistical literature 
\cite{Kel70,CHS81,FKN90,FZ90}
through their characteristic functions and their stochastic representations, which naturally follow from the spherically symmetric distributions after affine transformations. A first systematic treatment of
circular complex elliptically symmetric (C-CES) distributions was provided in the engineering literature \cite{KL86} and further 
fully studied in \cite{Ric06}, and in the tutorial paper \cite{OTKV12}. Then, the general complex representation of the RES distributions, called noncircular CES (NC-CES) distributions, was introduced in \cite{EK06}.

The aim of this chapter is twofold. At first, a short overview of RES, C-CES, and NC-CES distributions (as complex representations of the RES distributions) is introduced with the aim of providing a common background for the other chapters of this book. Secondly, the main definitions and properties of these distributions are listed and shortly discussed.

This chapter is organized as follows. 
Section \ref{sec:Definition of real elliptically symmetric distributions} defines the RES
distributions equivalently through their characteristic functions and to their stochastic representations. Particular attention is paid to the absolutely continuous case and to the subclass of compound Gaussian distributions. Section \ref{sec:Definition of the complex elliptically symmetric distributions}
defines the C-CES and NC-CES distributions as complex representations of the RES distributions.
Section \ref{sec:Basic properties} presents basic properties related to moments, affine transformations, marginal and conditional distributions, and summation stability.
Then some well-known instances of RES distributions are provided with their main properties in Section
\ref{sec:Example of elliptically symmetric distributions}.
The joint estimation of the symmetry center and scatter matrix is briefly discussed in Section
\ref{sec:Parameter estimation}
through the SM and SCM estimators, ML, $M$-estimators, and Tyler's $M$-estimators, 
asymptotic Gaussian distribution of scatter $M$-estimators,
and Slepian-Bangs formula. Finally, Section \ref{sec:conclusion} briefly concludes this chapter. 

The following notations are used throughout this chapter. 
Matrices
and vectors are represented by bold upper case and bold lower case
characters, respectively. Vectors are by default in column
orientation, while the superscripts $T$, $H$, $*$, and ${\#}$ stand for transpose,
conjugate transpose, conjugate, and Moore Penrose inverse,
respectively. 
$({\bf a})_k$ and $({\bf A})_{k,\ell}$ denote the $k$ and $(k,\ell)$-th element of the vector ${\bf a}$ and the matrix ${\bf A}$, respectively.
$\e(.)$, $|.|$, and $\tra(.)$ are the expectation, determinant, and trace operators, respectively. 
${\bf I}$ is the identity matrix with the appropriate dimension.
$\ve ({\bf A})$ denotes the
``vectorization'' operator that turns a matrix ${\bf A}$ into a vector by
stacking the columns of the matrix one below another
and ${\rm vecs}({\bf A})$ is the vector  that is obtained from $\ve ({\bf A})$  by eliminating all supradiagonal elements of ${\bf A}$.
These vectors are used
in conjunction with the Kronecker product ${\bf A}\otimes{\bf B}$
as the block matrix whose $(i,j)$ block element is $a_{i,j}{\bf B}$,
and with the commutation matrix ${\bf K}$ and the duplication matrix ${\bf D}$
of appropriate dimension such that
 $\ve({\bf C}^T)={\bf K}\ve({\bf C})$ and
 $\ve({\bf A})={\bf D}{\rm vecs}({\bf A})$ where ${\bf A}$ is symmetric. 
 ${\bf e}_k$ is the vector of appropriate dimension with 1 in the $k$th position and zeros elsewhere.
 The acronyms r.v., p.d.f. and c.d.f. for respectively random variable, probability density function and cumulative distribution function are used.
 Finally,
$\Gamma(u)\pardef \int_{0}^{\infty}t^{u-1}e^{-t}dt$ is the Gamma function
with $\Gamma(k)=(k-1)!$ for $k \in \mathbb{N}$, $B(k,\ell)$  denotes the Beta function with
$B(k,\ell)= \frac{\Gamma(k)\Gamma(\ell)}{\Gamma(k+\ell)}$ and
${\rm Gam}(k,\theta)$ is the Gamma distribution of scale $\theta$ with p.d.f.  
$p(x)=\frac{1}{\Gamma(k)\theta^k}x^{k-1}\exp(-x/\theta)$.
$x =_d y$, $x_n \rightarrow_d D$ and $x \sim D$ mean that the r.v. $x$ and $y$ have the same distribution,
the sequence of r.v. $x_n$ converges in distribution to $D$
and
$x$ follows the distribution $D$, respectively.
The subscripts $r$ and $c$ are used to refer to the real and complex data cases, respectively.
%
%
\section{Definition of the real elliptically symmetric distributions}
\label{sec:Definition of real elliptically symmetric distributions}
\subsection{Characteristic function}
\label{sec:Characteristic function RES}
After earlier works on this topic, RES distributions were formalized by \cite{Kel70} and further studied by \cite{CHS81,FKN90,FZ90}.
They were first defined as affine transformations of spherically distributed r.v.
Then, by the uniqueness theorem (see, e.g., \cite[pp. 346-351]{Sho00}, they were alternatively defined by their characteristic functions.
\begin{definition}
\label{def:fonction characterisque RES}
An r.v. ${\bf x} \in \mathbb{R}^m$ is said  to have a RES distribution if there exists a vector
$\boldsymbol \mu \in \mathbb{R}^m$, an $m \times m$ symmetric positive semi-definite matrix ${\bf \Sigma}$ of rank $k\le m$ and a function $\phi_{r}(.): \mathbb{R}^+\rightarrow  \mathbb{R}$
called {\it symmetry center}, {\it scatter matrix} and {\it characteristic generator}, respectively, such that the characteristic function of ${\bf x}$ is of the form
\begin{equation}
\label{eq:fonction characterisque RES}
\Phi_x({\bf t})
\pardef
\e[\exp(i {\bf t}^T{\bf x})]
=\exp(i {\bf t}^T{\boldsymbol \mu})
\phi_{r}({\bf t}^T{\bf \Sigma}{\bf t}),
\hspace{1cm}{\bf t} \in \mathbb{R}^m.
\end{equation}
\end{definition}
We shall write ${\bf x} \sim {\rm RES}_m(\boldsymbol \mu,{\bf \Sigma}, \phi_{r})$ and note that the couple $({\bf \Sigma}, \phi_{r}(.))$ does not uniquely identify the distribution of ${\bf x}$ because $(c^2{\bf \Sigma}, \phi_{r}(./c^2))$ gives the same distribution.  This scale ambiguity is easily avoided by restricting the function $\phi_{r}(.)$ in a suitable way
(e.g., by fixing a moment as it is explained in Section
\ref{sec:Moments}), or by putting a constraint on the scatter matrix ${\bf \Sigma}$ (e.g., $\tra({\bf \Sigma})=m$). Note that for $m=1$, these distributions
coincide with the class of one-dimensional symmetric
distributions w.r.t. the symmetry center. 
\subsection{Stochastic representation}
\label{sec:Stochastic representation RES}
%
Equivalently to the definition \eqref{eq:fonction characterisque RES}, the RES distributed r.v. ${\bf x}$ can be defined from an affine function
\begin{equation}
\label{eq:Affine definition}
{\bf x} \pardef{\boldsymbol \mu}+{\bf A}{\bf x}_s
\end{equation}
of a $k$-dimensional spherically distributed r.v. ${\bf x}_s$, where ${\bf A}\in \mathbb{R}^{m\times k}$ is any square root
(${\bf \Sigma}={\bf A}{\bf A}^T$) of the scatter matrix
${\bf \Sigma}$ of rank $k$, and thus full column rank. 
Such spherically distributions are defined equivalently in the following \cite[Chap.2]{FKN90}
\begin{definition}
\label{def:distribution spherique}
An r.v. ${\bf x}_s \in \mathbb{R}^k$ is spherically distributed i.f.f.
\begin{itemize}
\item
${\bf x}_s=_d {\cal O}{\bf x}_s$ for arbitrary real-valued $k$-dimensional orthonormal matrix 
${\cal O}$,
\item
There exists a function $\phi_r(.):$ $\mathbb{R}^+ \mapsto \mathbb{R}$, such that the characteristic function of ${\bf x}_s$ is given by
\begin{equation}
\label{eq:fonction caracteristique spherique}
\Phi_{x_s}({\bf t})
=\phi_r(\|{\bf t}\|^2), \ {\bf t} \in \mathbb{R}^k,
\end{equation}
\item
For every ${\bf h}\in \mathbb{R}^k$,\ ${\bf h}^T{\bf x}_s=_d\|{\bf h}\|x_{s_i}$ 
with ${\bf x}_s=(x_{s_1},.., x_{s_i},..,x_{s_k})^T$,
\item
If ${\bf x}_s$ is absolutely continuous w.r.t. Lebesgue measure on $\mathbb{R}^k$, 
there exists a function\footnote{both functions $\phi_r(.)$ and $g_r(.)$ are generally parameterized by the dimension $k$ and in practice by a finite-dimensional parameter (see examples in Section \ref{sec:Example of elliptically symmetric distributions}).} 
$g_r(.)$: $\mathbb{R}^+ \mapsto \mathbb{R}^+$
such that the p.d.f. of ${\bf x}_s$ w.r.t. this measure is of the form 
\begin{equation}
\label{eq:pdf spherique}
p({\bf x}_s)=g_r(\|{\bf x}_s\|^2),
\end{equation}
where 
\begin{equation}
\label{eq:normalization}
\delta_{r,k} \pardef \int_0^{\infty}t^{k/2-1}g_{r}(t)dt
=\frac{\Gamma(k/2)}{\pi^{k/2}},
\end{equation}
ensuring that 
$p({\bf x}_s)$ integrates to one.
\item
There exists a non-negative r.v. $\mathcal{Q}_{r,k}$, and ${\bf u}_{r,k}$ that are independent where ${\bf u}_{r,k}$  is uniformly distributed on the unit real $k$-sphere
(${\bf u}_{r,k} \sim U(\mathbb{R}S^k))$ such that
\begin{equation}
\label{eq:representation stochastique spherique}
{\bf x}_s
=_d\sqrt{\mathcal{Q}_{r,k}}{\bf u}_{r,k}.
\end{equation}
\end{itemize}
\end{definition}

Consequently from \eqref{eq:Affine definition}
and \eqref{eq:fonction caracteristique spherique},
we find the characteristic function 
\eqref{eq:fonction characterisque RES} 
which defines the RES distribution since
\[
\Phi_x({\bf t})=\exp(i {\bf t}^T{\boldsymbol \mu})\Phi_{x_s}({\bf A}^T{\bf t})
=\exp(i {\bf t}^T{\boldsymbol \mu})
\phi_{r}({\bf t}^T{\bf \Sigma}{\bf t}).
\]
From \eqref{eq:Affine definition} and 
\eqref{eq:representation stochastique spherique}, we obtain the following
\begin{theorem}
\label{th:Stochastic representation RES}
${\bf x}$ is ${\rm RES}_m(\boldsymbol \mu,{\bf \Sigma}, \phi_{r})$ distributed, i.f.f. it admits
the following stochastic full-rank representation
\begin{equation}
\label{eq:Stochastic representation RES}
{\bf x}
=_d{\boldsymbol \mu}+ \sqrt{\mathcal{Q}_{r,k}}{\bf A}{\bf u}_{r,k}
={\boldsymbol \mu}+ \mathcal{R}_{r,k}{\bf A}{\bf u}_{r,k}.
\end{equation}
\end{theorem}
The r.v. $\mathcal{Q}_{r,k}$ and ${\cal R}_{r,k}\pardef \sqrt{\mathcal{Q}_{r,k}}$ are the 2nd-order modular and modular (or generating) variates of the r.v. ${\bf x}$, respectively.
We note that there is a one-to-one mapping between the characteristic generator $\phi_{r}$ and the c.d.f. $F_{R_r}$ of ${\cal R}_{r,k}$ (called generating c.d.f.). Thus we can also write
${\bf x} \sim {\rm RES}_m(\boldsymbol \mu,{\bf \Sigma}, F_{R_r})$ and, equivalently to Definition \ref{def:fonction characterisque RES}, we retrieve the scale ambiguity in the couple 
$({\bf A},\mathcal{R}_{r,k})$ in
\eqref{eq:Stochastic representation RES}. 
Theorem 
\ref{th:Stochastic representation RES}  provides an obvious mechanism to generate r.v. 
${\bf x}\sim {\rm RES}_m(\boldsymbol \mu,{\bf \Sigma}, F_{R_r})$: it only involves generating
${\cal R}_{r,k}$ according to its c.d.f. $F_{R_r}$ and 
${\bf u}_{r,k}=\frac{{\bf n}_{r,k}}{\|{\bf n}_{r,k}\|}$ where ${\bf n}_{r,k}$ is $k$-dimensional zero-mean Gaussian distributed r. v. with covariance ${\bf I}$ (${\bf n}_{r,k}\sim \mathbb{R}N_k({\bf 0},{\bf I})$). Moreover, the following important property follows from Theorem 
\ref{th:Stochastic representation RES}:
\begin{equation}
\label{eq:forme quadratique}
\mathcal{Q}_{r,k}=_d({\bf x}-{\boldsymbol \mu})^T{\bf \Sigma}^{\#}({\bf x}-{\boldsymbol \mu}).
\end{equation}
%
\subsection{The absolutely continuous case}
\label{sec:The absolutely continuous case RES}
From \eqref{eq:Affine definition} the r.v. ${\bf x}$ is  absolutely continuous w.r.t. Lebesgue measure on $\mathbb{R}^m$, i.f.f. ${\bf x}_s$ is too.
From \eqref{eq:representation stochastique spherique}, it is immediate to verify that this condition is satisfied i.f.f. $\mathcal{Q}_{r,k}$ or 
$\mathcal{R}_{r,k}$ is absolutely continuous w.r.t. Lebesgue measure on
$\mathbb{R}^+$.
In this case, the p.d.f. of ${\bf x}$ is defined on the $k$-dimensional subspace of $\mathbb{R}^m$ spanned by the range space of ${\bf A}$.

In the particular case where $k=m$, i.e., rank$({\bf \Sigma})=m$, 
${\bf A}$ is a non singular $m \times m$ square matrix.
From
the one to one mapping \eqref{eq:Affine definition} between ${\bf x}$  and ${\bf x}_s$, and the p.d.f. \eqref{eq:pdf spherique},
the p.d.f. 
of ${\bf x}$ on $\mathbb{R}^m$ can be expressed as:
\begin{equation}
\label{eq:pdf RES}
p({\bf x})
=|{\bf \Sigma}|^{-1/2}g_{r}[({\bf x}-{\boldsymbol \mu})^T 
{\bf \Sigma}^{-1}({\bf x}-{\boldsymbol \mu})].
\end{equation}
We note that unlike the notation used in, e.g., \cite{OTKV12}, \eqref{eq:pdf RES} does not explicitly include
the usual p.d.f. normalizing constant.
In \eqref{eq:pdf RES} $g_{r}(.)$:
$\mathbb{R}^+ \mapsto \mathbb{R}^+$ is an arbitrary function, called {\it density generator} such that from 
\eqref{eq:pdf spherique} and 
\eqref{eq:normalization}
\begin{equation}
\label{eq:delta rmg}
\delta_{r,m}\pardef \int_0^{\infty}t^{m/2-1}g_{r}(t)dt
=\frac{\Gamma(m/2)}{\pi^{m/2}}.
\end{equation}
Clearly from \eqref{eq:pdf RES}, the couple $({\bf \Sigma},g_r(.))$
 does not uniquely identify the distribution of ${\bf x}$ because 
$(c^2{\bf \Sigma},c^{m}g_r(.\ c^2))$ gives the same distribution.

We adopt the notation ${\rm RES}_m({\boldsymbol \mu},{\bf \Sigma},g)$
instead of ${\rm RES}_m({\boldsymbol \mu},{\bf \Sigma},\phi)$.
The level sets of $p({\bf x})$ are a family of hyper ellipsoids in $\mathbb{R}^m$
symmetrically centered at $\boldsymbol \mu$, where shape and orientation are determined by ${\bf \Sigma}$. This justifies the  terminology of symmetrical elliptical distributions.
Furthermore, the stochastic representation \eqref{eq:Stochastic representation RES} reduces to
\begin{equation}
\label{eq:Stochastic representation RES rang plein}
{\bf x}
=_d{\boldsymbol \mu}+ \sqrt{\mathcal{Q}_{r,m}}{\bf \Sigma}^{1/2}{\bf u}_{r,m}
=
{\boldsymbol \mu}+ \mathcal{R}_{r,m}{\bf \Sigma}^{1/2}{\bf u}_{r,m}.
\end{equation}
Here too, note that $({\bf \Sigma},\mathcal{Q}_{r,m})$ and $(c^2{\bf \Sigma},c^{-2}\mathcal{Q}_{r,m})$ give the same distribution of ${\bf x}$.
\eqref{eq:Stochastic representation RES rang plein} implies that \eqref{eq:forme quadratique} simplifies to
\begin{equation}
\label{eq:Q}
\mathcal{Q}_{r,m}
=_d
({\bf x}-{\boldsymbol \mu})^T {\bf \Sigma}^{-1}({\bf x}-{\boldsymbol \mu}).
\end{equation}
From 
\eqref{eq:pdf RES},
\eqref{eq:Stochastic representation RES rang plein}
and \eqref{eq:Q},
the p.d.f. of
$\mathcal{Q}_{r,m}$ and 
$\mathcal{R}_{r,m}$ are respectively 
\begin{equation}
\label{eq:pdf R Q RES}
p(q)
=\delta^{-1}_{r,m}q^{m/2-1}g_{r}(q)
\ \ 
\mbox{and}
\ \ 
p(r)
=2\delta^{-1}_{r,m}r^{m-1}g_{r}(r^2).
\end{equation}

Finally, we note that the RES distributions do not necessarily possess a p.d.f. w.r.t. Lebesgue measure on $\mathbb{R}^m$ even when ${\bf \Sigma}$ is not singular. Such an example is the $U(\mathbb{R}S^m)$ distribution which belongs to ${\rm RES}_m({\bf 0},{\bf I}, \phi_r)$ distributions, where the explicit (but somewhat involved) form of $\phi_r$ can be found in \cite[the. 2.51]{FZ90}.
%
\subsection{The subclass of compound-Gaussian distributions}
\label{sec:The subclass of compound Gaussian distributions RES}
An important subclass of RES distributions are the compound-Gaussian (CG) distributions, whose circular complex 
representations (denoted C-CCG) have been widely used in the engineering literature, for example, for modeling radar clutter \cite{WTW06}. An r.v. having CG distributions with zero symmetry center is also referred to as {\it spherically invariant random vectors} (SIRV) in the engineering literature (see, e.g., in \cite{CL81,CLR95,RWO93,Yao73}) 
and as {\it scale mixtures of normal distributions}
in the statistics literature \cite{AM74}\cite{MVM2006}.
These distributions are defined by their stochastic representation.
\begin{definition}
\label{def:stochastic representation CG}
An r.v. ${\bf x} \in \mathbb{R}^m$ is said  to have a real CG distribution (RCG) if it admits the following representation
\begin{equation}
\label{eq:Stochastic representation RCG}
{\bf x}
=_d{\boldsymbol \mu}+ \sqrt{\tau_{r}} {\bf n},
\end{equation}
for some positive real r.v. $\tau_r$ with c.d.f. $F_{\tau}$ 
(not related neither to dimension $m$ nor to rank $k$), called the {\it texture} independent of 
${\bf n} \sim \mathbb{R}N_m({\bf 0},{\bf \Sigma})$, called the {\it speckle}. 
The r.v. $\sqrt{\tau_r}$ is often called {\it mixing variable} with mixing distribution in the statistical literature (see e.g., \cite{MVM2006}).
We write 
${\bf x} \sim {\rm RCG}_m({\boldsymbol \mu},{\bf \Sigma},F_{\tau})$ to denote this case.
\end{definition}

Note that \eqref{eq:Stochastic representation RCG} can be rewritten as
\begin{equation}
\label{eq:Stochastic representation RCG b}
{\bf x}
=_d{\boldsymbol \mu}+ \sqrt{\tau_{r}}{\bf A} {\bf n}_0,
\end{equation}
where ${\bf n}_0 \sim \mathbb{R}N_k({\bf 0},{\bf I})$.
Then by recalling that ${\bf n}_0=\|{\bf n}_0\| {\bf u}_{r,k}$ with
$s \pardef\|{\bf n}_0\|^2 \sim \chi^2_k={\rm Gam}(k/2,2)$ and
${\bf u}_{r,k}\sim U(\mathbb{R}S^k)$,
and where $s$ and ${\bf u}_{r,k}$ are independent.
It follows that
the stochastic representation \eqref{eq:Stochastic representation RCG} can also be written as
\begin{equation}
\label{eq:Stochastic representation RCG c}
{\bf x}
=_d{\boldsymbol \mu}+ \mathcal{R}_{r,k}{\bf A} {\bf u}_{r,k},
\end{equation}
where the modular variate $\mathcal{R}_{r,k} \pardef\sqrt{\tau_{r}s}$ and ${\bf u}_{r,k}$ are independent.
Consequently, the RCG distributions form a subclass of the RES distributions.
Furthermore, ${\bf x} \sim {\rm RES}_m({\boldsymbol \mu},{\bf \Sigma},\phi)$ belongs to the set of RCG distributions i.f.f. there exists an r.v. $\tau_r$ such that the 2nd-order modular variate $\mathcal{Q}_{r,k}=\mathcal{R}_{r,k}^2$ satisfies 
$\mathcal{Q}_{r,k}=\tau_{r}s$, i.e., 
$\mathcal{Q}_{r,k}$ is a scale mixture of the ${\rm Gam}(k/2,2)$ distribution.
This means that the conditional distribution of $\mathcal{Q}_{r,k}$ given $\tau_r=\tau$ is the 
${\rm Gam}(k/2,2\tau)$ distribution and thus the p.d.f of $\mathcal{Q}_{r,k}$ is
\begin{equation}
\label{eq:pdf Q RES=RCG}
p(q)
=\int_0^{\infty}
\frac{1}{\Gamma(k/2)(2\tau)^{k/2}}q^{k/2-1}\exp(-q/(2\tau))
dF_{\tau}(\tau).
\end{equation}

The characteristic function of  $\Phi_{x_s}({\bf t})$ of a RCG distributed r.v. ${\bf x}$ defined by
\eqref{eq:Stochastic representation RCG} is given straightforwardly by
\begin{eqnarray}
\nonumber
\Phi_{x_s}({\bf t})
&=&
\exp(i {\bf t}^T{\boldsymbol \mu})
\int_{0}^{\infty}\e(\exp(i {\bf t}^T\sqrt{\tau}{\bf n})/\tau)d F_{\tau}(\tau)
\\
\nonumber
&=&
\exp(i {\bf t}^T{\boldsymbol \mu})
\int_{0}^{\infty}\exp\left(-\frac{\tau}{2} {\bf t}^T{\bf \Sigma}{\bf t}\right)d F_{\tau}(\tau)
\\
\label{eq:cf RCG}
&=&
\exp(i {\bf t}^T{\boldsymbol \mu})\phi_r( {\bf t}^T{\bf \Sigma}{\bf t})
\end{eqnarray}
where the characteristic generator 
\begin{equation}
\label{eq: characteristic generator RCG}
\phi_r(u)=\int_{0}^{\infty}\exp\left(-\frac{\tau}{2}u\right)d F_{\tau}(\tau)
\end{equation}
does not depend on the dimension $m$ nor on the rank $k$,
unlike the RES distributions which are not RCG whose characteristic generator can depend on it.

In the particular case where ${\bf \Sigma}$ is not singular $(k=m)$, 
the conditional distribution of 
${\bf x}$ given $\tau_r=\tau$ is the $\mathbb{R}N_m({\boldsymbol \mu},\tau{\bf \Sigma})$ distribution from
\eqref{eq:Stochastic representation RCG}. 
Consequently, the distribution of ${\bf x}$ is always continuous w.r.t. Lebesgue measure on $\mathbb{R}^m$  and its p.d.f. is given by
\begin{equation}
\label{eq:pdf x RES=RCG}
p({\bf x})
=
(2\pi)^{-m/2}|{\bf \Sigma}|^{-1/2}
\int_0^{\infty}
\tau^{-m/2}\exp\left(
-\frac{1}{2\tau}({\bf x}-{\boldsymbol \mu})^T 
{\bf \Sigma}^{-1}({\bf x}-{\boldsymbol \mu})\right)
dF_{\tau}(\tau).
\end{equation}
Note that the p.d.f. \eqref{eq:pdf x RES=RCG} can always be written in the form 
\eqref{eq:pdf RES} with the density generator
\begin{equation}
\label{eq:denssity generator RES=RCG}
g_r(t)
=
(2\pi)^{-m/2}
\int_0^{\infty}
\tau^{-m/2}\exp\left(
-\frac{t}{2\tau}\right)
dF_{\tau}(\tau),
\end{equation}
and similarly to the RES distributions, we are faced with scale ambiguity where
$(c^2{\bf \Sigma}, c^{-2}\tau_r,c^m g_r(. c^2))$ gives the same R-CG distribution.
Note that \eqref{eq:denssity generator RES=RCG} reduces to the density generator \eqref{eq:gauss_cg_dg}
of the Gaussian distribution when $\tau_r$ is a degenerate r.v. putting all the probability at $\tau_r=1$.
Note also that the $\epsilon$-contaminated Gaussian distribution belongs to the class of CG distributions and is obtained when $\tau_r$ is a discrete r.v. with $P(\tau_r=a^2)=\epsilon$ and $P(\tau_r=1)=1-\epsilon$,
where $(a^2,\epsilon)$ are parameters that control the heaviness of the tails as compared to the Gaussian distribution.
\section{Definition of the complex elliptically symmetric distributions}
\label{sec:Definition of the complex elliptically symmetric distributions}
\subsection{Characteristic function}
\label{sec:Characteristic function NC-CES}
An r.v. ${\bf x} \in \mathbb{C}^m$ is said  to have a noncircular complex elliptically symmetric (NC-CES) distribution (also called generalized complex elliptical  in \cite{OK14} if
the r.v. $\bar{\bf x}\pardef ({\rm Re}({\bf x})^T,{\rm Im}({\bf x})^T)^T \in \mathbb{R}^{2m}$ is RES distributed.
Denote the symmetry center and the scatter matrix (of rank $k\le 2m$) of $\bar{\bf x}$ by $\bar{\boldsymbol \mu} \in \mathbb{R}^{2m}$ and
$\bar{\bf \Sigma}  \in \mathbb{R}^{2m \times 2m}$, respectively.
Using the one-to-one mapping $\bar{\bf x} \mapsto \widetilde{\bf x}\pardef ({\bf x}^T,{\bf x}^H)^T= \sqrt{2}{\bf M}\bar{\bf x}$ where
${\bf M}\pardef  
\footnotesize{\frac{1}{\sqrt{2}}\left(\begin{array}{cc}
{\bf I}& i{\bf I}\\
{\bf I} &  -i{\bf I}\\
\end{array}
\right)}$ is unitary, we obtain
$\bar{\bf t}^T\bar{\boldsymbol \mu}={\rm Re}({\bf t}^H{\boldsymbol \mu})$,
$\bar{\bf t}^T\bar{\bf \Sigma}\bar{\bf t}
=\frac{1}{2}\widetilde{\bf t}^H({\bf M}\bar{\bf \Sigma}{\bf M}^H)\widetilde{\bf t}
=\frac{1}{4}\widetilde{\bf t}^H\widetilde{\bf \Sigma}\widetilde{\bf t}$ with
$\bar{\bf t}\pardef ({\rm Re}({\bf t})^T,{\rm Im}({\bf t})^T)^T$,
$\bar{\boldsymbol \mu}\pardef ({\rm Re}({\boldsymbol \mu})^T,{\rm Im}({\boldsymbol \mu})^T)^T$,
$\widetilde{\bf t}\pardef ({\bf t}^T,{\bf t}^H)^T$ and
$\widetilde{\bf \Sigma}\pardef 2{\bf M}\bar{\bf \Sigma}{\bf M}^H$
of rank $k$ structured as
\begin{equation}
\label{eq:extended scatter matrix}
\widetilde{\bf \Sigma}
=
\left(\begin{array}{cc}
{\bf \Sigma}& {\bf \Omega}\\
{\bf \Omega}^* &  {\bf \Sigma}^* \\
\end{array}
\right),
\end{equation}
where ${\bf \Sigma}$ and ${\bf \Omega}$ defined from 
$2{\bf M}\bar{\bf \Sigma}{\bf M}^H$,
are positive semi-definite Hermitian and complex symmetric matrices, respectively.
Consequently, we obtain the following theorem by the definition \eqref{eq:fonction characterisque RES}
\begin{theorem}
\label{def:fonction characterisque NC-CES}
The characteristic function of an NC-CES distributed r.v. ${\bf x}\in \mathbb{C}^m$ is of the form
\begin{equation}
\label{eq:fonction characterisque NC-CES}
\Phi_x({\bf t})
\pardef
\e[\exp(i \bar{\bf t}^T\bar{\bf x})]
=\exp(i {\rm Re}({\bf t}^H{\boldsymbol \mu}))
\phi_{c}
\left(\frac{1}{2}\widetilde{\bf t}^H\widetilde{\bf \Sigma}\widetilde{\bf t}\right)
\hspace{1cm}{\bf t} \in \mathbb{C}^m,
\end{equation}
where
\begin{equation}
\label{eq:phi r c}
  \phi_{c}(u)\pardef \phi_{r}\left(\frac{1}{2}u\right)  
\end{equation}
is the characteristic generator, and where
${\boldsymbol \mu}\in \mathbb{C}^m$ and 
$\widetilde{\bf \Sigma}\in \mathbb{C}^{2m\times 2m}$ denote respectively the symmetric center and the {\it extended scatter matrix} of the NC-CES distributed r.v. 
${\bf x}$.
\end{theorem}

We shall write ${\bf x}\sim {\rm CES}_{m}({\boldsymbol \mu},{\bf \Sigma},{\bf \Omega},\phi_{c})$.
In the particular case where  ${\bf \Omega}={\bf 0}$, 
the rank of $\widetilde{\bf \Sigma}$ is even with 
${\rm rank}(\widetilde{\bf \Sigma})=2{\rm rank}({\bf \Sigma})=k$
and
${\bf x}$ is C-CES distributed\cite{KL86,Ric06,OTKV12}. The term \textit{circular} is often dropped in the current terminology used in signal processing where the distribution of a C-CES r.v. is usually indicated as ${\bf x}\sim {\rm CES}_{m}({\boldsymbol \mu},{\bf \Sigma},\phi_{c})$. Moreover, 
from \eqref{eq:fonction characterisque NC-CES}, we get:
\begin{theorem}
\label{def:fonction characterisque C-CES}
The characteristic function of  a C-CES distributed r.v. ${\bf x}\in \mathbb{C}^m$ is of the form
\begin{equation}
\label{eq:fonction characterisque C-CES}
\Phi_x({\bf t})
=\exp(i {\rm Re}({\bf t}^H{\boldsymbol \mu}))
\phi_{c}({\bf t}^H{\bf \Sigma}{\bf t})
\hspace{1cm}{\bf t} \in \mathbb{C}^m,
\end{equation}
where
${\boldsymbol \mu}\in \mathbb{C}^m$ and 
${\bf \Sigma}\in \mathbb{C}^{m\times m}$ denote the symmetric center, and the scatter matrix of the C-CES distributed r.v. ${\bf x}$, respectively.
\end{theorem}
%
\subsection{Stochastic representation}
\label{sec:Stochastic representation CES}
%
From the definition of the NC-CES distribution, a simple complex-valued extension of the stochastic representation \eqref{eq:Stochastic representation RES}
is only possible if the rank of $\widetilde{\bf \Sigma}$, which is equal to the rank of $\bar{\bf \Sigma}$, is even (it is, in particular, the case of the C-CES distribution and the case where $\widetilde{\bf \Sigma}$ is not singular for which $k=2m$).
Let $2k$ be the rank of $\widetilde{\bf \Sigma}$. In this case,
there exists an $m \times k$ full column rank matrix ${\bf A}$ such that ${\bf \Sigma}={\bf A}{\bf A}^H$
and ${\bf \Omega}={\bf A}{\bf \Delta}_{\kappa}{\bf A}^T$
where ${\bf \Delta}_{\kappa}=\diag(\kappa_1,\dots, \kappa_k)$ is a real diagonal matrix with non-negative real entries
$(\kappa_i)_{i=1,..,k}$ \cite[Corollary 4.6.12(b)]{HJ85}. Furthermore, it has been proved in \cite{EK06}, that
$0 \le \kappa_i \le 1$.
This parameterization allows us to state that the stochastic representation of this distribution, proved in \cite{AD19}, is a multivariate extension of the univariate generation of NC-CES r.v. presented in \cite[sec. IV.C]{OEK11}.
\begin{theorem}
\label{res:Stochastic representation NC-CES}
${\bf x}$ is ${\rm CES}_{m}({\bf \boldsymbol{\mu}},{\bf \Sigma},{\bf \Omega},\phi_{c})$ distributed i.f.f. it admits the following stochastic representation
\begin{equation}
\label{eq:stochastic representation NC-CES}
{\bf x}=_{d}{\bf \boldsymbol \mu}+\sqrt{\mathcal{Q}_{c,k}}{\bf A}[{\bf \Delta}_1{\bf u}_{c,k}+{\bf \Delta}_{2}{\bf u}^*_{c,k}],
\end{equation}
where $\mathcal{Q}_{c,k}\pardef \frac{1}{2}\mathcal{Q}_{r,2k}$ and ${\bf u}_{c,k}\sim U(\mathbb{C}S^k)$ are independent, 
${\bf \Delta}_1\pardef\frac{{\bf \Delta}_{+}+{\bf \Delta}_{-}}{2}$
and ${\bf \Delta}_2\pardef\frac{{\bf \Delta}_{+}-{\bf \Delta}_{-}}{2}$
where ${\bf \Delta}_{+} \pardef \sqrt{{\bf I}+{\bf \Delta}_{\kappa}}$
and ${\bf \Delta}_{-} \pardef \sqrt{{\bf I}-{\bf \Delta}_{\kappa}}$.
\end{theorem}

In the particular case of C-CES distributions, ${\bf \Omega}={\bf 0}$, which is equivalent to ${\bf \Delta}_{\kappa}={\bf 0}$, i.e., ${\bf \Delta}_{1}={\bf I}$ and 
${\bf \Delta}_{2}={\bf 0}$ and consequently the stochastic representation \eqref{eq:stochastic representation NC-CES} reduces to the well known
stochastic representation reported in \cite{OTKV12}:
\begin{equation}
\label{eq:Stochastic representation C-CES}
{\bf x}=_{d}{\bf \boldsymbol \mu}+\sqrt{\mathcal{Q}_{c,k}}{\bf A}{\bf u}_{c,k}.
\end{equation}

Note that similarly to the RES distribution, the C-CES distribution can be defined  from the affine function \eqref{eq:Affine definition}, where here 
${\bf x}_s$ is $k$-dimensional spherically distributed defined  by the equality ${\bf x}_s=_d {\cal U}{\bf x}_s$ for arbitrary complex-valued $k$-dimensional unitary matrix 
${\cal U}$, and where ${\bf A}\in \mathbb{C}^{m\times k}$ is any square root (${\bf \Sigma}={\bf A}{\bf A}^H$) of the scatter matrix
\cite{KL86},\cite{Ric06}.
Such r.v. ${\bf x}_s$  is also characterized by the stochastic representation
${\bf x}_s
=_d\sqrt{\mathcal{Q}_{c,k}}{\bf u}_{c,k}$
where the non-negative r.v. $\mathcal{Q}_{c,k}$, and ${\bf u}_{c,k}$ are independent, 
with ${\bf u}_{c,k} \sim U(\mathbb{C}S^k)$. ${\bf x}_s$ is also characterized by 
a characteristic function of the form 
$\Phi_{x_s}({\bf t})= \phi_{c}(\|{\bf t}\|^2)$, ${\bf t} \in \mathbb{C}^k$.
Consequently, from \eqref{eq:Affine definition}, we find the characteristic function 
\eqref{eq:fonction characterisque C-CES} which defines the C-CES distribution since
$\Phi_x({\bf t}) \pardef
\e[\exp(i {\rm Re}({\bf t}^H{\bf x}))]
=\exp(i {\rm Re}({\bf t}^H{\boldsymbol \mu}))\Phi_{x_s}({\bf A}^H{\bf t})
=\exp(i {\rm Re}({\bf t}^H{\boldsymbol \mu}))
\phi_{c}({\bf t}^H{\bf \Sigma}{\bf t})$.
\subsection{The absolutely continuous case}
\label{sec:The absolutely continuous case CES}
The p.d.f. of ${\bf x}$ is defined on $\mathbb{C}^{m}$ i.f.f. $\bar{\bf x}$ is absolutely continuous w.r.t. Lebesgue  measure on  $\mathbb{R}^{2m}$.
Assuming that rank($\widetilde{\bf \Sigma})=2m$ and
using the identities
$(\bar{\bf x}-\bar{\boldsymbol \mu})^T \bar{\bf \Sigma}^{-1}(\bar{\bf x}-\bar{\boldsymbol \mu})
=(\widetilde{\bf x}-\widetilde{\boldsymbol \mu})^H \widetilde{\bf \Sigma}^{-1}(\widetilde{\bf x}-\widetilde{\boldsymbol \mu})$ and
$| \bar{\bf \Sigma}|=2^{-2m}|\widetilde{\bf \Sigma}|$,
 the p.d.f. \eqref{eq:pdf RES} becomes
\begin{equation}
\label{eq:pdf NC CES}
p({\bf x})
=|\widetilde{\bf \Sigma}|^{-1/2}g_{c}
\left[\frac{1}{2}(\widetilde{\bf x}-\widetilde{\boldsymbol \mu})^H \widetilde{\bf \Sigma}^{-1}(\widetilde{\bf x}-\widetilde{\boldsymbol \mu})\right],
\end{equation}
where $g_{c}(t)$ is defined by
\begin{equation}
\label{eq:g r c}
g_{c}(t) \pardef 2^m g_{r}(2t),
\end{equation}
and $g_{r}(t)$ is the density generator associated with the distribution 
${\rm RES}_{2m}(\bar{\boldsymbol \mu},\bar{\bf \Sigma},\phi_{r})$,
which satisfies
\begin{equation}
\label{eq:normalization C}
\delta_{c,m}\pardef \int_0^{\infty}t^{m-1}g_{c}(t)dt =\delta_{r,2m} =\frac{\Gamma(m)}{\pi^m}. 
\end{equation}
We note that in this case, 
\eqref{eq:stochastic representation NC-CES} 
is written equivalently in the form
$\widetilde{\bf x}
=_d\widetilde{\boldsymbol \mu}
+\widetilde{\bf \Sigma}^{1/2}\widetilde{\bf u}_{c,m}$
(where $\widetilde{\bf u}_{c,m}\pardef({\bf u}_{c,m}^T,{\bf u}_{c,m}^H)^T$), then it follows that 
\begin{equation}
\label{eq:Q r c}
\mathcal{Q}_{c,m}=\frac{1}{2}\mathcal{Q}_{r,2m}
\end{equation}
and
\begin{equation}
\label{eq:Qcm}
\mathcal{Q}_{c,m}
=_d
\frac{1}{2}(\widetilde{\bf x}-\widetilde{\boldsymbol \mu})^H \widetilde{\bf \Sigma}^{-1}(\widetilde{\bf x}-\widetilde{\boldsymbol \mu}).  
\end{equation}

For C-CES distributed ${\bf x}$, \eqref{eq:pdf NC CES} reduces to
\begin{equation}
\label{eq:pdf C-CES}
p({\bf x})
=|{\bf \Sigma}|^{-1}g_{c}[({\bf x}-{\boldsymbol \mu})^H {\bf \Sigma}^{-1}({\bf x}-{\boldsymbol \mu})],
\end{equation}
and \eqref{eq:Stochastic representation C-CES} implies
\begin{equation}
\label{eq:Q C-CES}
\mathcal{Q}_{c,m}
=_d
({\bf x}-{\boldsymbol \mu})^H {\bf \Sigma}^{-1}({\bf x}-{\boldsymbol \mu}).
\end{equation}
Following the same derivation as for the RES distribution, we get
the following p.d.f. of $\mathcal{Q}_{c,m}$ and $\mathcal{R}_{c,m}$
\begin{equation}
\label{eq:pdf Q R C CES}
p(q)
=\delta^{-1}_{c,m}q^{m-1}g_{c}(q)
\ \ 
\mbox{and}
\ \
p(r)
=2\delta^{-1}_{c,m}r^{2m-1}g_{c}(r^2).
\end{equation}
\subsection{The subclass of compound Gaussian distributions}
\label{sec:The subclass of compound Gaussian distributions CES}
In many engineering applications, only the circular complex case is considered.
A r.v. ${\bf x}\in \mathbb{C}^m$ is said to be circular complex compound Gaussian (C-CCG) distributed
if the r.v. $\bar{\bf x}\in \mathbb{R}^{2m}$ is RCG distributed 
(see definition
\ref{def:stochastic representation CG}) where the associated extended scatter matrix
\eqref{eq:extended scatter matrix} is bloc-diagonal
$\widetilde{\bf \Sigma}
=\left(\begin{array}{cc}
{\bf \Sigma}& {\bf 0}\\
{\bf 0} &  {\bf \Sigma}^* \\
\end{array}\right)$.
Consequently from the one to one mapping $\bar{\bf x}\mapsto\tilde{\bf x}$,
definition \eqref{def:stochastic representation CG} gives
the following stochastic representation
\begin{equation}
\label{eq:Stochastic representation C-CCG}
{\bf x}
=_d{\boldsymbol \mu}+ \sqrt{\tau_{c}} {\bf n},
\end{equation}
where $\tau_{c}=2\tau_{r}$ 
with $\tau_{r}$, independent from ${\bf n}$, is associated with the $2m$-dimensional RCG distribution 
and ${\bf n} \sim \mathbb{C}N_m({\bf 0},{\bf \Sigma})$.

As a consequence, all the properties given in
Section \ref{sec:The subclass of compound Gaussian distributions RES} can be deduced from this real to complex representation. In particular with ${\bf \Sigma}={\bf A}{\bf A}^H$ of rank $k$:
${\bf n}_0 \sim \mathbb{C}N_k({\bf 0},{\bf I})$,
\eqref{eq:Stochastic representation RCG c} with now
$\mathcal{R}_{c,k} \pardef\sqrt{\tau_{c}s}$ where 
$s \pardef\|{\bf n}_0\|^2 \sim \frac{1}{2}\chi^2_{2k}={\rm Gam}(k,1)$.
Then, the C-CCG distributions form a subclass of the C-CES distributions and
${\bf x} \sim$ C-CES$_m({\boldsymbol \mu},{\bf \Sigma},\phi)$ belongs to the set of C-CCG distributions i.f.f. the p.d.f of $\mathcal{Q}_{c,k}$ is
\begin{equation}
\label{eq:pdf Q CES=CCG}
p(q)
=\int_0^{\infty}
\frac{1}{\Gamma(k){\tau}^k}q^{k-1}\exp(-q/{\tau})
dF_{\tau}(\tau)
\end{equation}
and \eqref{eq:pdf x RES=RCG} becomes
\begin{equation}
\label{eq:pdf x CES=CCG}
p({\bf x})
=
{\pi}^{-m}|{\bf \Sigma}|^{-1}
\int_0^{\infty}
\tau^{-m}\exp\left(
-\frac{1}{\tau}({\bf x}-{\boldsymbol \mu})^H
{\bf \Sigma}^{-1}({\bf x}-{\boldsymbol \mu})\right)
dF_{\tau}(\tau).
\end{equation}
%
\section{Basic properties}
\label{sec:Basic properties}
%
In this section and throughout the rest of this chapter, we mainly consider the RES distributions, knowing that the C-CES distributions are only a particular representation of them for even $m$, where the complementary scatter matrix ${\bf \Omega}$ defined in \eqref{eq:extended scatter matrix} is zero.
Consequently we drop the indices $r$ in $\delta_{r,k}$,
$g_{r}$, $\phi_r$, $\mathcal{Q}_{r,k}$, $\mathcal{R}_{r,k}$ and ${\bf u}_{r,k}$ associated with the v.a. ${\bf x}$.
We will show that these distributions benefit from most of the properties of the Gaussian distribution except the additive stability,
whose conditions are more restrictive (see e.g., the quick surveys in  \cite{MVM2003},\cite{Fra14}).
%
\subsection{Moments}
\label{sec:Moments}
%
From the full-rank stochastic representation
\eqref{eq:Stochastic representation RES}, it is clear that ${\bf x}$
admits $p$th-order moments i.f.f. $\e(\mathcal{R}_{k}^p)< \infty$.
Using the characteristic function \eqref{eq:fonction characterisque RES}, $\e(\mathcal{R}_{k}^p)< \infty$ i.f.f. 
the characteristic generator
$\phi({\bf t})$ is $p$ times differentiable.
In this case, the $p$th-order moments of ${\bf x}$ are given by 
$\e(x_1^{p_1}x_2^{p_2}...x_m^{p_m})
=\frac{1}{i^p}
\frac{{\partial^p \Psi}_x({\bf t})}
{\partial t_1^{p_1}\partial t_2^{p_2}...\partial t_m^{p_m}}|_{{\bf t}={\bf 0}}$ 
with $p=\sum_{i=1}^m p_i$ and ${\bf x}=(x_1,x_2,...,x_m)^T$.

Assuming that the correspondent moments are finite, one has
\begin{eqnarray}
\label{eq:mean}
\e({\bf x})
&=&
\boldsymbol \mu
\\
\label{eq:covariance}
\cov({\bf x})
&=&
\frac{\e(\mathcal{R}_{k}^2)}{k}{\bf \Sigma}
\ = \frac{\e(\mathcal{Q}_{k})}{k}{\bf \Sigma}
\ =-2\phi'({\bf 0}){\bf \Sigma}.
\end{eqnarray}
In particular for RCG distributions \eqref{eq:covariance} becomes
\begin{equation}
\label{eq:covariance CG} 
\cov({\bf x})
=
\e(\tau){\bf \Sigma}.
\end{equation}

We see from \eqref{eq:covariance} and \eqref{eq:covariance CG} that under finite second-order moment assumption,
the covariance of ${\bf x}$ does not necessarily coincide with the scatter matrix 
${\bf \Sigma}$, but these two matrices are proportional.
Note that many
second-order signal processing methodologies, such as, for example, subspace-based processing (where $k=m$), require an estimate of the covariance only at
up to a multiplicative scalar.
In this case, the {\it shape matrix}, defined as a scaled version of the scatter matrix
\begin{equation}
\label{eq:scale matrix} 
{\bf V}_s\pardef \frac{1}{s({\bf \Sigma)}}{\bf \Sigma}
\end{equation}
can be adopted to characterize the correlation structure. 
The scalar factor $s({\bf \Sigma)}$ must follow the conditions $s(a{\bf \Sigma)}=as({\bf \Sigma)}$, $\forall a>0$ and $s({\bf I)}=1$.
Even if the choice of the scale functionals $s({\bf \Sigma)}$ is entirely arbitrary, in signal processing literature, the most popular scale is the one on the trace of the scatter matrix, i.e., $s({\bf \Sigma)} \pardef \tra({\bf \Sigma})/m$ leading to the following shape matrix ${\bf V}_s= \frac{m}{\tra({\bf \Sigma})}{\bf \Sigma}$.
One can also find $s({\bf \Sigma)}\pardef[{\bf \Sigma}]_{1,1}$ and
$s({\bf \Sigma}) \pardef |{\bf \Sigma}|^{1/m}$, leading to $[{\bf V}_s]_{1,1}=1$ and $|{\bf V}_s|=1$, respectively.

Otherwise, \eqref{eq:covariance} and \eqref{eq:covariance CG} can be used to resolve the scale ambiguity of the couple 
$({\bf \Sigma}, \phi(.))$ in the definition \eqref{eq:fonction characterisque RES} of the RES distribution by fixing the constraint on the characteristic generator $\phi(.)$
\begin{equation}
\label{eq:contrainte sur Q}
\e(\mathcal{R}_{k}^2)=\e(\mathcal{Q}_{k})=k={\rm rank}({\bf \Sigma}) \ \ 
\mbox{or}\ \ \e(\tau)=1\ \ \mbox{for RCG distributions},
\end{equation}
which ensures that $\cov({\bf x})={\bf \Sigma}$.

When the r.v. ${\bf x}$ is absolutely continuous w.r.t.
Lebesgue measure on $\mathbb{R}^m$ and $k=m$, 
\eqref{eq:contrainte sur Q} is equivalent to
the following constraint on the density generator $g(.)$ thanks to \eqref{eq:pdf R Q RES}
\begin{equation}
\label{eq:contrainte sur g}
\delta_{m}^{-1} \int_0^{\infty}q^{m/2}g(q)dq
=m.
\end{equation}
If $\e(\mathcal{Q}_{m})$ is not finite, rather imposing 
Median$(\mathcal{R}_{m})=1$, i.e., from \eqref{eq:pdf R Q RES} the constraint 
$2\delta^{-1}_{m}\int_{0}^{1}r^{m-1}g(r^2)dr=\frac{1}{2}$ or equivalently
$\delta^{-1}_{m}\int_{0}^{1}t^{m/2-1}g(t)dt=\frac{1}{2}$,
is a more appropriate scaling constraint as its avoids any finite moment assumptions. 
Similarly for RCG distributions, if $\e(\tau)$ is not finite, the constraint 
Median$(\tau)=1$ (i.e., $F_{\tau}(1)=\frac{1}{2}$) can be used.
Indeed many RES distributions do not have finite second-order moments.

To consider higher-order multivariate central moments, let us consider $\sigma_{i_1,i_2,..,i_{\ell}}
\pardef \e[(x_{i_1}-\mu_{i_1})
(x_{i_2}-\mu_{i_2})...
(x_{i_{\ell}}-\mu_{i_{\ell}})]$ with 
${\boldsymbol \mu}= (\mu_1,\mu_2,...,\mu_m)^T$ and $(i_1,i_2,..,i_{\ell})\in \{1,..,m\}^{\ell}$.
By symmetry, all odd-order central moments are zero, provided that the corresponding moments do exist.
As for fourth-order moments, if $\e(\mathcal{R}_{m}^4)< \infty$ (with $k=m$), they all satisfy the identity \cite[p. 2]{Pai14}
\begin{equation}
\label{eq:ordre quatre}
\sigma_{i,j,k,\ell}
=(\kappa+1)
(\sigma_{i,j}\sigma_{k,\ell}
+\sigma_{i,k}\sigma_{j,\ell}
+\sigma_{i,\ell}\sigma_{j,k}), \ \ (i,j,k,\ell) \in \{1,...,m\}^4
\end{equation}
where 
\begin{equation}
\label{eq:kappa}
\kappa
\pardef
\frac{1}{3}
\left(
\frac{\sigma_{i,i,i,i}}{\sigma_{i,i}^2}-3
\right)
=
\frac{m}{m+2}\frac{\e(\mathcal{R}_{m}^4)}{(\e(\mathcal{R}_{m}^2))^2}-1
=\frac{\phi"(0)}{(\phi'(0))^2}-1.
\end{equation}
$\kappa$ is the {\it kurtosis} parameter of the marginal r.v. $x_i$ \cite{And84,Ben83}. 
It usually depends on the dimension $m$, but, remarkably, the kurtosis of the $i$th component does not depend on $i$, nor on the scatter matrix 
${\bf \Sigma}$. 
Consequently taken ${\bf \Sigma}={\bf I}$, we get for all RCG distributions: $\sigma_{i,i,i,i}=3\e(\tau^2)$ and $\sigma_{i,i}=\e(\tau)$ and thus
\begin{equation}
\label{eq:kappa RCG}
\kappa=\frac{\var(\tau)}{[\e(\tau)]^2}.
\end{equation}
Consequently the kurtosis parameter does not depend on the dimension $m$ for RCG distributions.
Note that $\kappa=0$ for the Gaussian distribution and that there exist other definitions of the kurtosis parameter or coefficient in the literature as
$\kappa \pardef  \frac{\sigma_{i,i,i,i}}{\sigma_{i,i}^2}-3$ and
$\kappa \pardef \frac{\sigma_{i,i,i,i}}{\sigma_{i,i}^2}$. 
Note also that the kurtosis is bounded below such that $\kappa \ge -2/(m+2)$ \cite{BB86}.
%
\subsection{Affine transformations and marginal distributions}
\label{sec:Affine transformations and marginal distributions}
%
Let ${\bf b} \in \mathbb{R}^n$ and 
${\bf B} \in \mathbb{R}^{n  \times m}$, and consider the transformed r.v. ${\bf y}={\bf B}{\bf x}+{\bf b}$. Its characteristic function 
$\Phi_y({\bf t})$ is deduced from the characteristic function \eqref{eq:fonction characterisque RES} of ${\bf x}$
 by
\begin{equation}
\label{eq:fonction characterisque RES de y}
\Phi_y({\bf t})
=
\exp(i {\bf t}^T{\bf b})
\Phi_x({\bf B}^T{\bf t})
=\exp(i {\bf t}^T({\bf B}{\boldsymbol \mu}+{\bf b}))
\phi({\bf t}^T{\bf B}{\bf \Sigma}{\bf B}^T{\bf t}),\ \ {\bf t}\in \mathbb{R}^n.
\end{equation}
Consequently, ${\bf y}$ is 
${\rm RES}_{n}({\bf B}{\boldsymbol \mu}+{\bf b},{\bf B}{\bf \Sigma}{\bf B}^T,\phi)$-distributed, and thus the class of 
${\rm RES}_{m}({\boldsymbol \mu},{\bf \Sigma},\phi)$ distributions is closed under affine transformations.
Note that the parameters are transformed as 
$({\boldsymbol \mu},{\bf \Sigma}) \mapsto 
({\bf B}{\boldsymbol \mu}+{\bf b},{\bf B}{\bf \Sigma}{\bf B}^T)$, which is the usual transformation for the couple (expectation, covariance) (when it exists) for the affine transformation 
${\bf x} \mapsto {\bf B}{\bf x}+{\bf b}$.
Note also that for $n \neq  m$, ${\bf y}$ may not belong to the same family as that of the r.v. ${\bf x}$ because
the characteristic generator $\phi$ may depend on the dimension $m$.

From the full-rank stochastic representation of ${\bf x}$
\eqref{eq:Stochastic representation RES}, we derive
\begin{equation}
\label{eq:Stochastic representation RES y}
{\bf y}
=_d {\bf B}{\boldsymbol \mu}+{\bf b} +\mathcal{R}_{k}({\bf B}{\bf A}){\bf u}_{k},
\end{equation}
which is not necessarily a full-rank stochastic representation of ${\bf y}$
because
${\rm rank}({\bf BA}) \le \min ({\rm rank}({\bf B}),{\rm rank}({\bf A})) \le \min (n, k)$.
If $\ell \pardef {\rm rank}({\bf BA})$, there exist full column rank $n \times \ell$ matrices ${\bf C}$
such that ${\bf B}{\bf \Sigma}{\bf B}^T={\bf C}{\bf C}^T$. 
And thus, similarly  to ${\bf x}$, where its distribution is equivalently defined by its
 stochastic full rank representation \eqref{eq:Stochastic representation RES} and by its characteristic function \eqref{eq:fonction characterisque RES}, we get from 
\eqref{eq:fonction characterisque RES de y}
a stochastic full rank representation of ${\bf y}$
\begin{equation}
\label{eq:Stochastic representation RES y full rank}
{\bf y}
=_d{\bf B}{\boldsymbol \mu}+{\bf b} + \mathcal{R}_{\ell}{\bf C}{\bf u}_{\ell},
\end{equation}
where the non-negative r.v. $\mathcal{R}_{\ell}$, and ${\bf u}_{\ell}$ are independent, 
${\bf u}_{\ell}$  is uniformly distributed on the unit real $\ell$-sphere.

Of course, it directly follows that univariate and multivariate marginals of
${\bf x}$ also are 
RES distributed with $\phi$ remains unchanged. 
If ${\bf x}=(x_1,...,x_m)^T=({\bf x}_1^T,{\bf x}_2^T)^T$, 
${\boldsymbol \mu}=(\mu_1,...,\mu_m)^T=({\boldsymbol \mu}_1^T,{\boldsymbol \mu}_2^T)^T$
and ${\bf \Sigma}=({\sigma}_{i,j})_{i,j=1}^m
=\footnotesize{\left(\begin{array}{cc}
{\bf \Sigma}_{1,1}& {\bf \Sigma}_{1,2}\\
{\bf \Sigma}_{2,1} &  {\bf \Sigma}_{2,2}\\
\end{array}
\right)}$, where ${\bf x}_1$ and ${\boldsymbol \mu}_1$ are $m_1$-dimensional vectors, 
${\bf \Sigma}_{1,1}$ and ${\bf \Sigma}_{2,2}$ are $m_1 \times m_1$ and $m_2 \times m_2$ matrices, respectively (with $m_1+m_2=m$), then
${\bf x}_1 \sim {\rm RES}_{m_1}({\boldsymbol \mu}_1,{\bf \Sigma}_{1,1},\phi)$,
${\bf x}_2 \sim {\rm RES}_{m_2}({\boldsymbol \mu}_2,{\bf \Sigma}_{2,2},\phi)$,
and
$x_i \sim {\rm RES}_{1}(\mu_i,\sigma_{i,i},\phi)$, $i=1,...,m$.
Of course, their stochastic representations \eqref{eq:Stochastic representation RES y} and 
\eqref{eq:Stochastic representation RES y full rank}
also, follow. For example for $m_1 \ge k$, if ${\bf A}
=
\footnotesize{\left[\begin{array}{c}
{\bf A}_1\\
{\bf A}_2\\
\end{array}\right]}$ 
where ${\bf A}_1$ is a  $m_1 \times k$ full column rank matrix
\begin{equation}
\label{eq:Stochastic representation RES x1 full rank}
{\bf x}_1
=_d
{\boldsymbol \mu}_1+ \mathcal{R}_{k}{\bf A}_1{\bf u}_{k}
\end{equation}
is a  stochastic full-rank representation of the marginal ${\bf x}_1$.
As for the univariate marginal $x_i$, the stochastic representations 
\eqref{eq:Stochastic representation RES y} and 
\eqref{eq:Stochastic representation RES y full rank}
reduce to 
\begin{equation}
\label{eq:Stochastic representation RES xi}
x_i
=_d
\mu_i+ \mathcal{R}_{k}{\bf a}_i^T{\bf u}_{k}
=_d
\mu_i+ \mathcal{R}_{1}\|{\bf a}_i\|u_{1},
\end{equation}
where ${\bf a}_i$ is the $i$-th column of ${\bf A}^T$ 
(thus $\|{\bf a}_i\|^2=\sigma_{i,i}$),
$\mathcal{R}_{k}$ and $\mathcal{R}_{1}$ are the modular variates of ${\bf x}$
and $x_i$, respectively,
and where $u_{1}$ reduces to the uniform discrete r.v. $\{-1,+1\}$.

Now we take a closer look at the marginal distributions when $k=m$.
In this case, 
the stochastic full-rank representation 
\eqref{eq:Stochastic representation RES x1 full rank}
and \eqref{eq:Stochastic representation RES xi}
of arbitrary univariate or multivariate marginal r.v.
${\bf x}_1$ and $x_i$ reduce to 
\begin{equation}
\label{eq:Stochastic representation RES k=m}
{\bf x}_1
=_d
{\boldsymbol \mu}_1+ \mathcal{R}_{m_1}{\bf \Sigma}_{1,1}^{1/2}{\bf u}_{m_1}
\ \ 
\mbox{and}
\ \
x_i
=_d
\mu_i+ \mathcal{R}_{1}\sigma_{i,i}^{1/2}u_{1},
\end{equation}
where the modular variates $\mathcal{R}_{m_1}$ of ${\bf x}_1$
(which includes the modular variates of the univariate marginal r.v. $x_i$ for $m_1=1$) and $\mathcal{R}_{m}$ of ${\bf x}$ are related by the relation \cite[Corollary p.59]{FZ90}
\begin{equation}
\label{eq:relation modular variates for k=m}
\mathcal{R}_{m_1}
=_d
\mathcal{R}_{m}\
\beta_{\frac{m_1}{2},\frac{m_2}{2}}.
\end{equation}
In \eqref{eq:relation modular variates for k=m} the r.v. $\mathcal{R}_{m}$ and
$\beta_{\frac{m_1}{2},\frac{m_2}{2}}$ are independent and $\beta_{\frac{m_1}{2},\frac{m_2}{2}}^2\sim {\rm Beta}(\frac{m_1}{2},\frac{m_2}{2})$.

Moreover in the absolutely continuous case w.r.t. Lebesgue measure on $\mathbb{R}^m$, i.e., 
${\bf x} \sim {\rm RES}_m({\boldsymbol \mu},{\bf \Sigma},g)$, 
\eqref{eq:relation modular variates for k=m} allows us to relate the p.d.f. $p_{m_1}(r)$ of 
$\mathcal{R}_{m_1}$ to the p.d.f. $p_{m}(r)$ of $\mathcal{R}_{m}$ \cite[rel. 2.5.15]{FZ90}
\begin{equation}
\label{eq:pdf R m1 et R m}
p_{m_1}(r)
=\frac{2 r^{m_1-1}}{{\rm B}(\frac{m_1}{2},\frac{m_2}{2})}
\int_{r}^{+\infty}t^{-(m-2)}(t^2-r^2)^{\frac{m_2}{2}-1}p_m(t)dt
\end{equation}
and to the density generator $g(.)$ of ${\bf x}$ using $p_{m}(r)$ given by \eqref{eq:pdf R Q RES}
\begin{equation}
\label{eq:pdf R m1 et g}
p_{m_1}(r)
=\frac{2 \pi^{m/2} r^{m_1-1}}{{\Gamma}(\frac{m_1}{2}){\Gamma}(\frac{m_2}{2})}
\int_{r^2}^{+\infty}(t-r^2)^{\frac{m_2}{2}-1}g(t)dt.
\end{equation}

This allows us to derive the density generators $g_{m_1|m}(.)$
 of the multivariate and univariate marginal r.v. ${\bf x}_1$ and $x_i$
thanks to \eqref{eq:pdf R Q RES} $p_{m_1}(r)=2\delta_{m_1}^{-1}r^{m_1-1}g_{m_1|m}(r^2)$
\begin{equation}
\label{eq:pdf m1 et g}
g_{m_1|m}(u)
=\delta_{m_2}^{-1}
\int_{u}^{+\infty}(t-u)^{\frac{m_2}{2}-1}g(t)dt.
\end{equation}

Therefore the p.d.f. of the multivariate marginal r.v. ${\bf x}_{1} \sim {\rm RES}_{m_1}({\boldsymbol \mu}_1,{\bf \Sigma}_{1,1},g_{m_1|m})$
and
$x_i \sim {\rm RES}_1(\mu_i,\sigma_{i,i},g_{1|m})$ 
are given respectively by
\begin{eqnarray}
\label{eq:pdf x m1}
p_{m_1|m}({\bf x}_1)
&=&
|{\bf \Sigma}_{1,1}|^{-1/2}g_{m_1|m}[({\bf x}_1-{\boldsymbol \mu}_1)^T 
{\bf \Sigma}_{1,1}^{-1}({\bf x}_1-{\boldsymbol \mu}_1)].
\\
\label{eq:pdf x i}
p_{i|m}(x_i)
&=&
\frac{1}{\sqrt{\sigma_{i,i}}}g_{1|m}\left(\frac{(x_i-\mu_i)^2}{\sigma_{i,i}}\right).
\end{eqnarray}

Note that there is no guarantee that the integral in \eqref{eq:pdf m1 et g} leads to a closed-form expression when a closed-form expression of $g(u)$ is available, except for the class of RCG distributions for which $g_{m_1|m}(u) = g_{m_1}(u)$, where $g_{m_1}(u)$ denotes the density generator $g(u)$ of ${\bf x}$ evaluated at $m=m_1$. In particular for the Gaussian distribution, we obtain  $g_{m_1|m}(u)=\frac{1}{(2\pi)^{m_1/2}}\exp(-\frac{u}{2})$.

For the RCG distributions defined in Subsection
\ref{sec:The subclass of compound Gaussian distributions RES}, it is clear from the stochastic representation 
\eqref{eq:Stochastic representation RCG} that 
\begin{equation}
\label{eq:Stochastic representation CG xi}
{\bf x}_1
=_d
{\boldsymbol \mu}_1 + \sqrt{\tau} {\bf n}_1,
\end{equation}
with ${\bf n}_1 \sim \mathbb{R}N_{m_1}({\bf 0},{\bf \Sigma}_{1,1})$.
Consequently, all the marginals belong to the same subclass of RCG distributions with the same c.d.f $F_{\tau}(.)$
and the density generator \eqref{eq:pdf m1 et g} reduces thanks to \eqref{eq:denssity generator RES=RCG} to
\begin{equation}
\label{eq:marginal denssity generator RES=RCG}
g_{m_1|m}(t)
=
(2\pi)^{-m_1/2}
\int_0^{\infty}
\tau^{-m_1/2}\exp\left(
-\frac{t}{2\tau}\right)
dF_{\tau}(\tau).
\end{equation}
In fact this property characterizes the  RCG distributions. This point is specified by a consistency result \cite[Th. 1]{Kan94}. 
This result states that any univariate and multivariate marginal distribution of an r.v. ${\bf x}$
belong to the same family as that of the  r.v. ${\bf x}$, i.f.f. the RES distribution belongs to the class of RCG distributions. 
This condition is also equivalent to the characteristic generator $\phi$ not related to the dimension $m$.
This consistency result shows that not all elliptically symmetric distributions can be used to define random processes.
Indeed from Kolmogorov's theorem, only among the elliptically symmetric distributions, the $m$ dimensional RGG distributed r.v. ${\bf x}=(x_1,...,x_m)^T$ for all $m$ can define a unique random process $(x_n)_{n \in  \mathbb{N}}$ that is characterized by the distribution of $\tau$, symmetric center and scatter of the process \cite{DA2024}.
%
\subsection{Conditional distributions}
\label{sec:Conditional distributions}
%
Let ${\bf x}=({\bf x}_1^T,{\bf x}_2^T)^T \sim 
{\rm RES}_{m}({\boldsymbol \mu},{\bf \Sigma},\phi)$ where ${\bf x}_1 \in \mathbb{R}^{m_1}$.
The conditional distribution of ${\bf x}_2$  given ${\bf x}_1$ is generally more difficult to describe than its marginal distribution presented in Section \ref{sec:Affine transformations and marginal distributions}.
For the sake of simplicity, we consider here only the full rank case ($k=m$);
see e.g. \cite{CHS81} and \cite[Th. 7, Cor. 8]{Fra14} for more general statements. 
In this case, it is proved \cite[Th 2.18]{FKN90}, that ${\bf x}_2$  given ${\bf x}_1={\bf x}_1^0$ (denoted 
${\bf x}_2|{\bf x}_1={\bf x}_1^0$)  is 
${\rm RES}_{m_2}({\boldsymbol \mu}_{2|1},{\bf \Sigma}_{2|1},\phi_{2|1})
 ={\rm RES}_{m_2}({\boldsymbol \mu}_{2|1},{\bf \Sigma}_{2|1},F_{2|1})$
 distributed with
 \begin{eqnarray}
\label{eq:esperance conditonnel}
{\boldsymbol \mu}_{2|1}
&=&
{\boldsymbol \mu}_{2}
+{\bf \Sigma}_{2,1}{\bf \Sigma}_{1,1}^{-1}
({\bf x}_1^0-{\boldsymbol \mu}_{1})
\\
\label{eq:scatter conditonnel}
{\bf \Sigma}_{2|1}
&=&
{\bf \Sigma}_{2,2}-{\bf \Sigma}_{2,1}
{\bf \Sigma}_{1,1}^{-1}
{\bf \Sigma}_{1,2},
\end{eqnarray}
and where $\phi_{2|1}$ and $F_{2|1}$ correspond respectively to the characteristic generator of ${\bf x}_2|{\bf x}_1={\bf x}_1^0$ 
and the c.d.f. of the conditional modular variate 
${\cal R}_{2|1}$ defined by the following stochastic full rank representation
\begin{equation}
\label{eq: conditional stochastic full rank representation} 
({\bf x}_2|{\bf x}_1={\bf x}_1^0)
=_d
{\boldsymbol \mu}_{2|1}+
\mathcal{R}_{2|1}{\bf A}_{2|1}{\bf u}_{m_2},
\end{equation}
where ${\bf A}_{2|1}$ is an arbitrary square root of ${\bf \Sigma}_{2|1}$ (i.e.,
${\bf \Sigma}_{2|1}={\bf A}_{2|1}{\bf A}_{2|1}^T$ and 
$\mathcal{R}_{2|1}$ and ${\bf u}_{m_2}$ are independent.
$\mathcal{R}_{2|1}$ is given by
\begin{equation}
\label{eq: conditional modular} 
\mathcal{R}_{2|1}
=_d
[\mathcal{R}_m^2-({\bf x}_1^0-{\boldsymbol \mu}_1)^T{\bf \Sigma}_{1,1}^{-1}({\bf x}_1^0-{\boldsymbol \mu}_1)]^{1/2}|{\bf x}_1={\bf x}_1^0.
\end{equation}

Otherwise if the conditional covariance $\cov({\bf x}_2|{\bf x}_1={\bf x}_1^0)$ exists, its expressions can be derived from \eqref{eq:covariance} and we get
\begin{equation}
\label{eq: conditional covariance} 
\cov({\bf x}_2|{\bf x}_1={\bf x}_1^0)
=\frac{1}{m_2}\left(\e[
\mathcal{R}_m^2|{\bf x}_1={\bf x}_1^0]-({\bf x}_1^0-{\boldsymbol \mu}_1)^T{\bf \Sigma}_{1,1}^{-1}({\bf x}_1^0-{\boldsymbol \mu}_1)\right)
{\bf \Sigma}_{2|1}.
\end{equation}

We note that the expressions of the conditional symmetry center 
${\boldsymbol \mu}_{2|1}$ and 
conditional scatter matrix ${\bf \Sigma}_{2|1}$
 are those obtained for the Gaussian distribution, but the conditional characteristic generator $\phi_{2|1}$ no longer belongs to the same family of RES, except for the Gaussian distribution.
 
For RCG distributions, the p.d.f. of ${\bf x}_2|{\bf x}_1={\bf x}_1^0$, (denoted $p({\bf x}_2/{\bf x}_1^0)$, is given by
\begin{equation}
\label{eq: ddp conditionnelle 0} 
p({\bf x}_2/{\bf x}_1^0)=\frac{1}{p({\bf x}_1^0)}\int_{0}^{\infty}p({\bf x}_2/{\bf x}_1^0,\tau)p({\bf x}_1^0/\tau) dF_{\tau}(\tau)
\end{equation}
where
$p({\bf x}_1^0/\tau)
=(2\pi \tau)^{-m_1/2}|{\bf \Sigma}_{1,1}|^{-1/2}\exp(-\frac{1}{2\tau}d^2_{{\bf \Sigma}_{1,1}}({\bf x}_1^0,{\boldsymbol \mu}_1))$
and
$p({\bf x}_2/{\bf x}_1^0,\tau)$

\noindent
$=(2\pi \tau)^{-m_2/2}|{\bf \Sigma}_{2|1}|^{-1/2}\exp(-\frac{1}{2\tau}d^2_{{\bf \Sigma}_{2|1}}({\bf x}_2,{\boldsymbol \mu}_{2|1}))$
with 
$d^2_{{\bf \Sigma}_{1,1}}({\bf x}_1^0,{\boldsymbol \mu}_1)
\pardef ({\bf x}_1^0-{\boldsymbol \mu}_1)^T{\bf \Sigma}_{1,1}^{-1}({\bf x}_1^0-{\boldsymbol \mu}_1)$
and
$d^2_{{\bf \Sigma}_{2|1}}({\bf x}_2,{\boldsymbol \mu}_{2|1})
\pardef ({\bf x}_2-{\boldsymbol \mu}_2)^T{\bf \Sigma}_{2|1}^{-1}({\bf x}_2-{\boldsymbol \mu}_2)$.
This yields
\begin{eqnarray}
\nonumber
p({\bf x}_2/{\bf x}_1^0)
&=&
\frac{(2\pi)^{-m_2/2}|{\bf \Sigma}_{2|1}|^{-1/2}}{\int_{0}^{\infty}\tau^{-m_1/2}
\exp(-\frac{1}{2\tau}d^2_{{\bf \Sigma}_{1,1}}({\bf x}_1^0,{\boldsymbol \mu}_1)) dF_{\tau}(\tau)}
\\
\label{eq: ddp conditionnelle 1} 
&&
\hspace{-2.5cm}
\times
\int_{0}^{\infty}\tau^{-m/2}
\exp\left(-\frac{1}{2\tau}d^2_{{\bf \Sigma}_{2|1}}({\bf x}_2,{\boldsymbol \mu}_{2|1}\right)
\exp\left(-\frac{1}{2\tau}d^2_{{\bf \Sigma}_{1,1}}({\bf x}_1^0,{\boldsymbol \mu}_1)\right) dF_{\tau}(\tau).
\end{eqnarray}
Comparing \eqref{eq: ddp conditionnelle 1} to \eqref{eq:pdf RES}, 
we check that ${\bf x}_2|{\bf x}_1={\bf x}_1^0$ is 
${\rm RES}_{m_2}({\boldsymbol \mu}_{2|1},{\bf \Sigma}_{2|1},F_{2|1})$ distributed.
Moreover, comparing \eqref{eq: ddp conditionnelle 1} to \eqref{eq:pdf x RES=RCG}, we see
that ${\bf x}_2|{\bf x}_1={\bf x}_1^0$ is  RCG distributed,
but with a differently distributed texture $\tau$ than ${\bf x}$ and the marginal ${\bf x}_1$.
%
\subsection{Summation stability}
\label{sec:Summation stability}
%
Consider now the sum ${\bf y}$ of $n$ independent r.v. 
${\bf x}_1,..,{\bf x}_i,..,{\bf x}_n$ with the same scatter matrix, where 
${\bf x}_i \sim {\rm RES}_{m}({\boldsymbol \mu}_i,{\bf \Sigma},\phi_i)$. 
The characteristic function $\Phi_y({\bf t})$ of 
${\bf y}=\sum_{i=1}^n {\bf x}_i$ is
\begin{equation}
\label{eq: fonction characteristique somme} 
\Phi_y({\bf t})
=\prod_{i=1}^n \Phi_{x_i}({\bf t})
=
\exp\left(i {\bf t}^T(\sum_{i=1}^n {\boldsymbol \mu}_i)\right)
\prod_{i=1}^n\phi_i({\bf t}^T{\bf \Sigma}{\bf t})
=
\exp(i {\bf t}^T{\boldsymbol \mu})
\phi({\bf t}^T{\bf \Sigma}{\bf t}),
\end{equation}
where
${\boldsymbol \mu}
\pardef{\sum_{i=1}^n {\boldsymbol \mu}_i}$
and 
$\phi(u)
\pardef 
\prod_{i=1}^n\phi_i(u)$, which is structured as \eqref{eq:fonction characterisque RES}.
Consequently, the sum ${\bf y}$ is RES-distributed too \cite{Fra14}.
Similarly, for independent univariate r.v. $x_1,...,x_i,...,x_n$ but with arbitrary scatters $\sigma_i$, where
$x_i \sim {\rm RES}_{1}(\mu_i,\sigma_{i},\phi_i)$, 
the characteristic function $\Phi_y(t)$ of 
$y=\sum_{i=1}^n x_i$ is given by
$\Phi_y(t)
=
\exp(i t \mu)
\phi(t^2)$
where 
$\mu
\pardef \sum_{i=1}^n \mu_i$
and 
$\phi(u)
\pardef 
\prod_{i=1}^n\phi_i(\sigma_i u)$ and then 
$y \sim {\rm RES}_{1}(\mu,1,\phi)$. The sum $y$ is then symmetrically distributed w.r.t. $\mu$.

However, it is worth underlying that if the r.v. 
${\bf x}_1,..,{\bf x}_i,..,{\bf x}_n$ belong to the same family of RES distribution
(i.e., with $\phi_i=\phi$, $i=1,..,n$), 
the sum is not of the same family except for the so-called elliptical $\alpha$-stable distributions
\cite{NOL13}, (i.e., of characteristic functions
$\Phi_x({\bf t})
=
\exp(i {\bf t}^T{\boldsymbol \mu})
\exp(-\frac{1}{2}({\bf t}^T{\bf \Sigma}{\bf t})^{\alpha/2})$ with $\alpha \in (0,2]$), which includes the Gaussian distribution for $\alpha=2$.
Finally note that if the scatter matrices ${\bf \Sigma}_i$ of 
${\bf x}_1,..,{\bf x}_i,..,{\bf x}_n$ are not identical, the summation stability is generally lost for independent multivariate non-Gaussian, RES distributed r.v..

The independence condition of the v.a. ${\bf x}_i$ can be relaxed by the following properties proved in \cite[th. 4.2]{HL02}.
If the full-rank stochastic representations
\eqref{eq:Stochastic representation RES}
${\bf x}_i
=_d{\boldsymbol \mu}_i+ \mathcal{R}_{k}^i{\bf A}{\bf u}_{k}^{i}$, $i=1,2$ with 
${\bf \Sigma}={\bf A}{\bf A}^T$ satisfy the condition
$(\mathcal{R}_{k}^1,\mathcal{R}_{k}^2)$, 
${\bf u}_{k}^{1},{\bf u}_{k}^{2}$ are mutually independent, whereas 
$\mathcal{R}_{k}^1,\mathcal{R}_{k}^2$ 
may be dependent on each other, the sum ${\bf y}$ is also 
${\rm RES}_{m}({\boldsymbol \mu},{\bf \Sigma},\phi)$-distributed, where the expression of $\phi$ is given in \cite[th. 4.2]{HL02}. $\phi$ reduces to the product $\phi_1\phi_2$ when $\mathcal{R}_{k}^1$ and $\mathcal{R}_{k}^2$ are independent. A natural application of this property is in the context of a multivariate time series.
%
%
\section{Example of elliptically symmetric distributions}
\label{sec:Example of elliptically symmetric distributions}
%
In this section, we present some examples of elliptically symmetric distributions and discuss their main specific properties. Throughout this section, we mainly consider the case of distributions  absolutely continuous w.r.t. Lebesgue  measure on  $\mathbb{R}^{m}$ with rank$({\bf \Sigma})=m$.
Each distribution is defined by its density generator under a functional form parameterized by $m$ and a finite-dimensional parameter.
We mainly use their real-valued definition through the p.d.f. \eqref{eq:pdf RES}
where the characteristic and density generators, and the texture \eqref{eq:Stochastic representation RCG} are simply denoted here by respectively $\phi$, $g$ and $\tau$, knowing that C and NC complex-valued definitions \eqref{eq:pdf NC CES} and \eqref{eq:pdf C-CES} are simply deduced from the real-valued definition with double dimension
\eqref{eq:phi r c}, \eqref{eq:g r c}. For C-CES distributions, interested readers can consult \cite{OTKV12}.
%
\subsection{Gaussian distribution}
\label{sec:Gaussian distribution}
The Gaussian distribution is the best-known and widely used distribution among the RES class in classical signal processing applications. Its ubiquity is mainly due to the central limit theorem (CLT) that allows for the use of the Gaussian distribution as a good and handy approximation of the statistical behavior of a set of observations in many practical scenarios. 
		
As a particular case of RES distribution, the Gaussian distribution 
denoted hereafter by $\mathbb{R}N_m(\boldsymbol \mu,{\bf \Sigma})$,
is characterized by a \textit{characteristic} and a \textit{density} generators that can be expressed respectively as:
\begin{equation}
\label{eq:gauss_cg_dg}
\phi(u)=\exp \left( -\frac{u}{2} \right) 
\ \ 
\mbox{and}
\ \ 
g(t)=\frac{1}{(2\pi)^{m/2}}\exp \left( -\frac{t}{2} \right).
\end{equation}
Using \eqref{eq:pdf R Q RES}, when ${\bf \Sigma}$ is non-singular, it is immediately verified that the p.d.f. of the second-order modular variate $\mathcal{Q}$ of a Gaussian r.v. is given by:
\begin{equation}
\label{eq:pdf_Q_gauss}
p(q)=\frac{\delta^{-1}_{m}}{(2\pi)^{m/2}}q^{m/2-1}\exp \left( -\frac{q}{2}\right)  = \frac{1}{2^{m/2}\Gamma(m/2)}q^{m/2-1}\exp \left( -\frac{q}{2}\right), \; q \in \mathbb{R}^+,
\end{equation}
where $\delta_{m} = \frac{\Gamma(m/2)}{\pi^{m/2}}$ from \eqref{eq:normalization}. 
It is immediate to verify that the p.d.f. in \eqref{eq:pdf_Q_gauss} is the one of a central $\chi^2$-distribution with $m$ degrees of freedom, i.e., $\mathcal{Q} \sim \chi^2_m$. Note that this result is perfectly in line with the well-known property of the Gaussian r.v. with mean value $\boldsymbol{\mu}$ and covariance matrix $\boldsymbol{\Sigma}$ whose quadratic form $(\mathbf{x}-\boldsymbol{\mu})^T \boldsymbol{\Sigma}^{-1}(\mathbf{x}-\boldsymbol{\mu})$ is $\chi^2_m$-distributed.
	
Using again \eqref{eq:pdf R Q RES}, we can express the p.d.f. of the modular variate $\mathcal{R} \sim \sqrt{\chi^2_m}$ as :
\begin{equation}
\label{eq:pdf_R_gauss}
p(r)=\frac{1}{2^{m/2-1}\Gamma(m/2)}r^{m-1}\exp \left( -\frac{r^2}{2}\right), \; r \in \mathbb{R}^+.
\end{equation}

Note that from the real to complex representation given in Section 
\ref{sec:Definition of the complex elliptically symmetric distributions},
the density generator of complex Gaussian distributions becomes from \eqref{eq:g r c}
\begin{equation}
\label{eq:density generator gauss complex}
g(t)=\frac{1}{{\pi}^{m}}\exp (-t)
\end{equation}
and consequently, the p.d.f. of the circular and non-circular Gaussian distribution, denoted respectively 
by $\mathbb{C}N_m(\boldsymbol \mu,{\bf \Sigma})$ and $\mathbb{C}N_m(\boldsymbol \mu,{\bf \Sigma},{\bf \Omega})$
are given from \eqref{eq:pdf C-CES} and \eqref{eq:pdf NC CES} by respectively
\begin{eqnarray}
\nonumber
p({\bf x})
&=&
\pi^{-m}|{\bf \Sigma}|^{-1}\exp[({\bf x}-{\boldsymbol \mu})^H {\bf \Sigma}^{-1}({\bf x}-{\boldsymbol \mu})],
\\
\label{eq:pdf complex gaussian}
p({\bf x})
&=&
\pi^{-m}|\widetilde{\bf \Sigma}|^{-1/2}\exp
\left[\frac{1}{2}(\widetilde{\bf x}-\widetilde{\boldsymbol \mu})^H \widetilde{\bf \Sigma}^{-1}(\widetilde{\bf x}-\widetilde{\boldsymbol \mu})\right],
\end{eqnarray}
where ${\bf \Sigma}$ and ${\bf \Omega}$ are the covariance and complementary covariance (or pseudo covariance) of ${\bf x}$, respectively.

The Gaussian distribution is usually used as a reference to define heavier-tailed and lighter-tailed distributions in the RES class. Some examples will be given below. 
%
\subsection{Student's $t$-distribution}
\label{sec:Student's t  distribution}
A popular example of heavy-tailed distribution is the Student $t$-distribution (or simply $t$-distribution). The $t$-distribution belongs to the RES class, and it is characterized by the following density generator \cite{robust_t_dist}: 
\begin{equation}
\label{eq:t_dg}
g(u) = 
\frac{\Gamma(\frac{\nu+m}{2})}{(\nu\pi)^{m/2}\Gamma({\nu}/2 )}\left( 1 + \frac{u}{\nu} \right)^{-\frac{\nu+m}{2}},
\; u\in \mathbb{R}^+ 
\end{equation}
with $0 < \nu < \infty$ \textit{degrees of freedom}. The parameter $\nu$ controls the tails of the distribution that are uniformly heavier than the Gaussian ones. In particular, for small values of $\nu$, the sampled r.v. are highly non-Gaussian while, as $\nu \rightarrow \infty$, the $t$-distribution collapses into the Gaussian one.
The case $\nu=1$ is called the Cauchy distribution.
		
As shown in \cite{robust_t_dist}, the second-order modular variate of an $m$-dimensional $t$-distributed r.v. satisfies the following property:
\begin{equation}
\label{eq:t_Q_F}
\mathcal{Q} / m \sim F_{m,\nu},
\end{equation}
where $F_{m,\nu}$ denotes the $F$-distribution with $m$ and $\nu$ degrees of freedom \cite{JKB1994}. Moreover, using \eqref{eq:pdf R Q RES}, the p.d.f. of the second-order modular variate $\mathcal{Q}$ and of the modular variate 
$\mathcal{R}$ can be directly derived.

It is worth mentioning that the $t$-distribution belongs to the subclass of the CG distribution described in Section \ref{sec:The subclass of compound Gaussian distributions CES}. Specifically, $t$-distributed r.v. can be obtained by generating the texture $\tau$ such that its inverse is distributed as a Gamma random variable,
\begin{equation}
\label{eq:tau t}{\tau}^{-1} \sim  {\rm Gam}(\nu/2, 2/\nu),
\end{equation}
Since the $t$-distribution belongs to the subclass of the CG distributions, it holds that the marginals belong to the same family.
Hence in particular the univariate marginals $x_i$ are $t$-distributed with parameter $\nu$, symmetry center 
$({\boldsymbol \mu})_i$ 
and scatter $({\bf \Sigma})_{i,i}$ with density generator given by \eqref{eq:t_dg} where $m=1$.

This distribution has first, second and fourth-order moments if, respectively, $\nu>1$, $\nu>2$, and $\nu>4$.
In particular, we have $\e(\tau) = \frac{\nu}{2(\nu-2)}$, $\e(\mathcal{Q})=\frac{m\nu}{\nu-2}$ 
and $\cov({\bf x})=\frac{\nu}{\nu-2}{\bf \Sigma}$ for $\nu>2$, while the kurtosis, defined in subsection \ref{sec:Moments} in \eqref{eq:kappa} is $\kappa = \frac{2}{\nu-4}$ for $\nu>4$.
%
\subsection{Generalized Gaussian distributions}
\label{sec:Generalized Gaussian distributions}
The Generalized Gaussian (GG) distribution is a RES distribution with the remarkable ability to characterize both heavy-tailed and light-tailed (with respect to the Gaussian one) data behavior. 
The multivariate GG distribution was introduced in \cite{MVM1998} as a multivariate generalization of the power exponential distribution.
As a member of the RES class, using the parametrization \cite[rel. (27)]{OTKV12} and \eqref{eq:g r c}, the GG distribution has a density generator given by:
\begin{equation}
\label{eq:GG_dg}
g(u) = 
\frac{s\Gamma(m/2)}{(2\pi)^{m/2}b^{m/2s}\Gamma({m/2s)}} \exp \left( -\frac{u^s}{2^sb}\right), \; u\in \mathbb{R}^+
\end{equation}
where $s > 0$ and $b > 0$ are two parameters generally called \textit{exponent} (or \textit{shape}) and \textit{scale}, respectively. 
The scale parameter $b$ is used to ensure that $p({\bf x})$ integrates to 1 in \eqref{eq:pdf RES}, while the exponent parameter $s$ controls the non-Gaussianity. In particular, for $0<s<1$, the tails of the GG distribution are heavier with respect to the Gaussian one, while for $s > 1$ they are lighter.
For $s$ tending to $\infty$, this distribution converges to a uniform distribution in an ellipsoid centered on ${\boldsymbol \mu}$.
Clearly, for $s=1$ we get the Gaussian distribution. 
Otherwise for $s=1/2$, we get the Laplace distribution (also called double exponential distribution). Note that \eqref{eq:GG_dg} is consistent
with \cite{Kotz1968} and \cite{PBTB2013} with different definitions of the scale.

It was proved in \cite{MVM2008} that for $ s\in (0,1]$, the GG distributions 
is a scale mixture of normal distribution, i.e., whose p.d.f. is written in the form \eqref{eq:pdf x RES=RCG}, where the p.d.f. of the mixing variable $\sqrt{\tau}$ is given by \cite[(2.2)]{MVM2008} and depends on $m$. So it does not belong to the subset of CG distributions  in the sense given by Definition
\ref{def:stochastic representation CG}.
For $s >1$, the GG distributions does not belong to the CG family either.
In particular, in this case, the univariate marginals of a multivariate GG distribution are not power exponential distributed, i.e., with p.d.f. given by \cite{Lindsey1999}
\begin{equation}
\label{eq:pdf power exponential}
p(x)
=\frac{s}{\sigma \sqrt{\pi} b^{1/2s}\Gamma(1/2s)}\exp\left(-\frac{1}{2b}\left(\frac{x-\mu}{\sigma}\right)^{2s}\right).
\end{equation}
with $b \pardef [\frac{1}{2}\Gamma(\frac{1}{2s})/\Gamma(\frac{3}{2s})]^s$ with $\sigma^2= \var(x)$.

Regarding the second-order modular variate $\mathcal{Q}$ of an $m$-dimensional GG-distributed r.v., it is straightforward to deduce from \eqref{eq:pdf R Q RES}
and \eqref{eq:GG_dg} that:
\begin{equation}
\label{eq:GG_Q_G}
\mathcal{Q}^s \sim {\rm Gam}\left( \frac{m}{2s},2^sb\right),
\end{equation}
We can exploit again \eqref{eq:pdf R Q RES} to get a closed-form expression of the p.d.f. of $\mathcal{Q}$ and  $\mathcal{R}$.

We note that the moments exist at all orders and $\cov({\bf x})={\bf \Sigma}$ if the parameters $s$ and $b$ are related by the relation $b=[\frac{m}{2}\Gamma(\frac{m}{2s})/\Gamma(\frac{m/2+1}{s})]^s$. 
 Note that in this case
the kurtosis $\kappa$ defined in 
\eqref{eq:kappa} depends on the dimension $m$, e.g., for $m=1$ we get 
$\kappa=\frac{\Gamma(5/2s)\Gamma(1/2s)}{3[\Gamma(3/2s)]^2}-1$ which gives $\kappa=0$ for $s=1$ (Gaussian case).
%
\subsection{The $K$-distribution}
\label{sec:The $K$-distribution}
Another example of RES distribution belonging to the CG-subclass is the $K$-distribution. An $m$-dimensional r.v. $\mathbf{x}$ is said to be $K$-distributed if it has the CG-representation \eqref{eq:Stochastic representation RCG}
characterized by a Gamma-distributed texture,
\begin{equation}
\label{eq:tau K}
\tau \sim \frac{1}{2}{\rm Gam}(\nu, 1/\nu),
\end{equation}
where $\nu > 0$ is a \textit{shape} parameter. So the p.d.f. of $\tau$ is 
\begin{equation}
\label{eq:pdf K tau}
p(\tau)=\frac{2\nu^{\nu}}{\Gamma(\nu)}(2\tau)^{\nu-1}\exp(-2\nu\tau),\ \tau\in \mathbb{R}^+
\end{equation}
with 
$\e(\tau)=1/2$.

By using the definition of CG-distribution, it can be shown that the density generator is given by:
\begin{equation}
\label{eq:K_dg}
g(u) = \frac{\nu^{m/2}}{2^{\nu - 1} \pi^{m/2}\Gamma(\nu)}(2\nu u)^{(2\nu-m)/4}K_{(2\nu-m)/2}\left( \sqrt{2\nu u}\right) ,\; u\in \mathbb{R}^+ ,
\end{equation} 
where $K_\alpha(\cdot)$ denotes the modified Bessel function of the second kind of order $\alpha$. 
The shape $\nu > 0$ is a parameter that controls the tails of the $K$-distribution: when $\nu \rightarrow 0$, the tails become heavier  while for $\nu \rightarrow \infty$ the $K$-distribution collapses onto the Gaussian one. Moreover, using \eqref{eq:pdf R Q RES}, the p.d.f. of the second-order modular variate $\mathcal{Q}$ and of the modular variate 
$\mathcal{R}$ can be directly derived.

Since the $K$-distribution belongs to the subclass of the CG distributions, it holds that the marginals belong to the same family.
Hence, in particular the univariate marginals $x_i$ are $K$-distributed with parameter $\nu$, 
symmetry center 
$({\boldsymbol \mu})_i$ 
and scatter $({\bf \Sigma})_{i,i}$
with density generator given by \eqref{eq:K_dg} where $m=1$.

To conclude, we note that the moments of all orders of the $K$-distribution exist and it can be shown that the kurtosis in \eqref{eq:kappa} is given by $\kappa = \frac{1}{\nu}$.
%
\subsection{Related distribution: the angular central Gaussian distribution}
\label{sec:Related distribution: the angular central Gaussian distribution}
To conclude this section, let us introduce a distribution that has a strong link with the RES family, even if it does not belong to this class.
An r.v. ${\bf x}$ is said angular central Gaussian (ACG) distributed if it admits the stochastic representation:
\begin{equation}
\label{eq:ACG gaussien}
\mathbf{x}  =_d \frac{\mathbf{n}}{\left\|\mathbf{n} \right\| }, \;
\mbox{where}\; \mathbf{n} \sim \mathbb{R}N_m(\boldsymbol{0}, \boldsymbol{\Sigma}).
\end{equation}
For non-singular $\boldsymbol{\Sigma}$, the p.d.f. of this distribution is given by 
\cite{Tyler_ACG}
\begin{equation}
\label{eq:ACG_pdf}
p(\mathbf{x}) =
 \frac{2\pi^{m/2}}{\Gamma(m/2)} \left| \boldsymbol{\Sigma} \right|^{-1/2}
 \left( \mathbf{x}^T \boldsymbol{\Sigma}^{-1}  \mathbf{x}\right)^{-m/2},
   \;  \mathbf{x} \in \mathbb{R}S^m.
\end{equation}
Even if this expression may look similar to the p.d.f. of a centered RES-distributed r.v., it is worth underlying a crucial difference: the density of a RES r.v., given in \eqref{eq:pdf RES}, is defined w.r.t. the Lebesgue measure on $\mathbb{R}^m$, while the density of an ACG r.v. in \eqref{eq:ACG_pdf} is defined w.r.t. the Lebesgue measure on real unit sphere $\mathbb{R}S^m$. For this reason, the stochastic representation provided in Sec. \ref{sec:Stochastic representation RES} does not hold for ACG r.v..
It can be noted from \eqref{eq:ACG gaussien} or \eqref{eq:ACG_pdf}, that the parameter
$\boldsymbol{\Sigma}$ can only be identified up to a multiplicative scalar factor.

Furthermore, it turns out that if ${\bf x}$ is arbitrarily centered RES distributed (not necessarily Gaussian), i.e., 
if $\mathbf{x} \sim \mathrm{RES}_m(\boldsymbol{0}, \boldsymbol{\Sigma}, \phi)$,
we get from its stochastic representation \eqref{eq:Stochastic representation RES} with 
$\boldsymbol{\Sigma}={\bf A}{\bf A}^T$ of rank $k$ and where 
${\bf u}_k=_d\frac{{\bf n}_0}{\left\|{\bf n}_0 \right\|}$ and
${\bf n}_0 \sim \mathbb{R}N_k({\bf 0},{\bf I})$
\begin{equation}
\label{eq:ACG RES}
\frac{{\bf x}}{\left\|{\bf x} \right\|}
=_d
\frac{{\bf A}{\bf u}_k}{\left\|{\bf A}{\bf u}_k \right\|}
=_d
\frac{{\bf A}{\bf n}_0}{\left\|{\bf A}{\bf n}_0 \right\|}
=_d
\frac{{\bf n}}{\left\|{\bf n} \right\| }, \;
\mbox{where}\; {\bf n} \sim \mathbb{R}N_m({\boldsymbol 0}, {\boldsymbol \Sigma}).
\end{equation}
Consequently the projection of ${\bf x}$ onto the unit real $m$-sphere is also 
ACG-distributed. In its definition \eqref{eq:ACG gaussien}, 
$\mathbb{R}N_m(\boldsymbol{0}, \boldsymbol{\Sigma})$ can be replaced by any 
$\mathrm{RES}_m(\boldsymbol{0}, \boldsymbol{\Sigma}, \phi)$ distribution
and the term ACG appears to be a slight misnomer.

Finally, note that this class of distribution is closed under standardized linear transformations, i.e., ${\bf x}\sim$ ACG$_m({\bf 0},\boldsymbol{\Sigma})$, then 
${\bf w}\pardef{\bf B}{\bf x}/\|{\bf B}{\bf x}\|\sim$ ACG$_k({\bf 0},{\bf B}\boldsymbol{\Sigma}{\bf B}^T)$ for any nonzero $k \times m$ matrix {\bf B}.
%
\section{Parameter estimation}
\label{sec:Parameter estimation}
%
We are interested in this Section to the estimation of the symmetry center ${\boldsymbol \mu}$ and scatter matrix ${\bf \Sigma}$ (assumed invertible here) of RES distributions.
Suppose we have an i.i.d. sample ${\bf x}_1,..,{\bf x}_i,..,{\bf x}_n$ of size $n>m$ from 
a ${\rm RES}_{m}({\boldsymbol \mu},{\bf \Sigma},g)$ distribution absolutely continuous w.r.t. Lebesgue measure on 
$\mathbb{R}^m$ under finite second-order moments.
The identifiability issue is solved here by imposing constraint 
\eqref{eq:contrainte sur Q} which ensures that
$\cov({\bf x}_i)={\bf \Sigma}$.
%
\subsection{Sample mean and sample covariance matrix}
\label{sec:Sample mean and sample covariance matrix}
%
A natural estimate of the 
parameters ${\boldsymbol \mu}$ and ${\bf \Sigma}$ are the sample mean 
$\widehat{\boldsymbol \mu}\pardef \frac{1}{n}\sum_{i=1}^n{\bf x}_i$ and the sample covariance matrix (SCM)
$\widehat{\bf \Sigma}\pardef
\frac{1}{n-1}\sum_{i=1}^n({\bf x}_i-\widehat{\boldsymbol \mu})({\bf x}_i-\widehat{\boldsymbol \mu})^T$, respectively.
It is well known that $\widehat{\boldsymbol \mu}$ and $\widehat{\bf \Sigma}$
are unbiased and mutually uncorrelated estimators. Under the particular case of Gaussian data, the two estimators are then independent.
$\widehat{\boldsymbol \mu}$ is RES distributed with symmetry center ${\boldsymbol \mu}$ and scatter matrix
$\frac{1}{n}{\bf \Sigma}$ by imposing constraint 
\eqref{eq:contrainte sur Q} on this distribution which does not necessarily belong to the same RES distribution family as ${\bf x}_i$ (see Section \ref{sec:Summation stability}).
Under finite fourth-order moments, applying the CLT to $\sum_{i=1}^n{\bf x}_i{\bf x}^T_i$ and $\sum_{i=1}^n{\bf x}_i$, it can be shown that $\widehat{\bf \Sigma}$ is asymptotically Gaussian distributed, i.e.,
\begin{equation}
\label{eq:th CLT Scatter}
\sqrt{n}\left(
\ve(\widehat{\bf \Sigma})-\ve({\bf \Sigma})
\right) \rightarrow_d \mathbb{R}N_{m^2}\left({\bf 0},{\bf R}_{{\Sigma}_{SCM}}  \right)
\end{equation}
with \cite[p. 5]{Pai14}
\begin{equation}
\label{eq:asymptotic covariance SCM}
{\bf R}_{{\Sigma}_{SCM}}
=(1+\kappa)({\bf I}+{\bf K})({\bf \Sigma}\otimes{\bf \Sigma})+\kappa \ve({\bf \Sigma})\ve^T({\bf \Sigma}),
\end{equation}
where $\kappa$ is the kurtosis parameter \eqref{eq:kappa} of the ${\rm RES}_{m}({\boldsymbol \mu},{\bf \Sigma},g)$ distribution. Note that this asymptotic distribution 
would also be the asymptotic distribution of the estimate 
$\frac{1}{n}\sum_{i=1}^n({\bf x}_i-{\boldsymbol \mu})({\bf x}_i-{\boldsymbol \mu})^T$ of ${\bf \Sigma}$ when
${\boldsymbol \mu}$ is known. 
For heavier tails than Gaussian distributions, $\kappa>0$ and can be very large without any upper-bound 
(see e.g. in Section \ref{sec:Student's t  distribution} when $\kappa$ approaches 4 for Student's $t$-distribution) and consequently, the SCM estimate can be a very bad estimator.
%
\subsection{ML estimation}
\label{sec:ML estimation}
%
To take into account the particular RES distribution of the data, the maximum likelihood (ML) estimator is often considered as the reference estimator because it is generally (when it exists and is unique) asymptotically efficient with a speed of convergence in $\sqrt{n}$. From Slepian-Bangs formula \eqref{eq:structure FIM}, the Fisher information matrix (FIM) for the parameter 
$({\boldsymbol \mu},{\rm vecs}({\bf \Sigma}))$ is given by
\begin{equation}
\label{eq:FIM}
{\rm FIM}
\left(\begin{array}{c}
{\boldsymbol \mu}\\
{\rm vecs}({\bf \Sigma})\\
\end{array}\right)
\!=\!
n\left(\begin{array}{cc}
a_0{\bf \Sigma}^{-1}& {\bf 0}\\\
{\bf 0} & {\bf D}^T[a_1({\bf \Sigma}^{-1}\otimes {\bf \Sigma}^{-1})\!+\!a_2\ve({\bf \Sigma}^{-1})\ve^T({\bf \Sigma}^{-1})]{\bf D}\\
\end{array}\right),
\end{equation}
where
\begin{equation}
\label{eq:a1a2a3}
a_0=\frac{\e[\mathcal{Q}\varphi^2(\mathcal{Q})]}{m}, \ \
a_1=\frac{\e[\mathcal{Q}^2\varphi^2(\mathcal{Q})]}{2m(m+2)} \ \ \mbox{and}\ \
a_2=\frac{1}{4}\left(\frac{\e[\mathcal{Q}^2\varphi^2(\mathcal{Q})]}{m(m+2)}-1\right),
\end{equation}
assuming that $g$ is continuously differentiable with $\varphi(t) \pardef -2g'(t)/g(t)$. 
Hence under existence, uniqueness and usual regularity conditions, the ML estimate
$(\widehat{\boldsymbol \mu},\widehat{\bf \Sigma})$ is
asymptotically Gaussian distributed:
\begin{equation}
\label{eq:th ML asymptotic}
\sqrt{n}
\left(\begin{array}{c}
\widehat{\boldsymbol \mu}-{\boldsymbol \mu}\\
\ve(\widehat{\bf \Sigma})-\ve({\bf \Sigma})
\end{array}
\right)
 \rightarrow_d \mathbb{R}N_{m+m^2}
\left(
\left(\begin{array}{c}
 {\bf 0} \\
{\bf 0} 
\end{array}\right),
\left(
\begin{array}{cc}
{\bf R}_{{\mu}_{ML}}& {\bf 0}\\
{\bf 0} & {\bf R}_{{\Sigma}_{ML}}\\
\end{array}\right)
\right),
\end{equation}
where ${\bf R}_{{\mu}_{ML}}$ and ${\bf R}_{{\Sigma}_{ML}}$ can be deduced from 
\eqref{eq:FIM} using the efficiency of the ML estimate:
\begin{equation}
\label{eq:Covariance ML}
{\bf R}_{{\mu}_{ML}}=\sigma_0{\bf \Sigma}
\ \ \mbox{and}\ \ 
{\bf R}_{{\Sigma}_{ML}}
=\sigma_1({\bf I}+{\bf K})({\bf \Sigma}\otimes {\bf \Sigma})
+\sigma_2\ve({\bf \Sigma})\ve^T({\bf \Sigma}),
\end{equation}
with \cite{tyler1982radial}
\begin{equation}
\label{eq:sigma012}
\sigma_0=\frac{m}{\e[\mathcal{Q}\varphi^2(\mathcal{Q})]},\
\sigma_1=\frac{m(m+2)}{\e[\mathcal{Q}^2\varphi^2(\mathcal{Q})]}
\ \ \mbox{and}\ \
\sigma_2=-\frac{2\sigma_1(1-\sigma_1)}{2+m(1-\sigma_1)},
\end{equation}
where $\sigma_1$ and $\sigma_2$ are free of scale ambiguity in contrast to $\sigma_0$.
Note that this asymptotic distribution was first given in \cite{tyler1982radial} by using the general structure of the covariance of random matrices whose distributions are radial.

Using the Delta method (see, e.g., \cite[chap. 6]{Sho00}), the asymptotic distribution of the ML estimate $\widehat{\bf V}_s$ of any shape matrix ${\bf V}_s$ defined by
\eqref{eq:scale matrix} from the scale $s({\bf \Sigma})$ can be deduced:
\begin{equation}
\label{eq:shape ML asymptotic}
\sqrt{n}
\left(
\ve(\widehat{\bf V}_s)-\ve({\bf V}_s)
\right)
 \rightarrow_d \mathbb{R}N_{m^2}
\left(
{\bf 0},
{\bf R}_{V_{s,ML}}
\right),
\end{equation}
with 
\begin{equation}
\label{eq:Covariance shape ML}
{\bf R}_{V_{s,ML}}
=\sigma_1
{\bf P}_s({\bf V}_s)
({\bf I}+{\bf K})({\bf V}_s\otimes {\bf V}_s){\bf P}_s^T({\bf V}_s),
\end{equation}
where 
${\bf P}_s({\bf V}_s)\pardef {\bf I}-{\rm vec}({\bf V}_s)\frac{d s({\bf \Sigma})}{d{\rm vec}^T({\bf \Sigma})}$ 
is given for 
 $s({\bf \Sigma})=[{\bf \Sigma}]_{1,1}$,
  $s({\bf \Sigma})=\frac{1}{m}\tra({\bf \Sigma})$ and 
  $s({\bf \Sigma})=|{\bf \Sigma}|^{1/m}$ by 
  ${\bf P}_s({\bf V}_s)={\bf I}-{\rm vec}({\bf V}_s){\bf e}^T_1$, 
${\bf P}_s({\bf V}_s)={\bf I}-\frac{1}{m}{\rm vec}({\bf V}_s){\rm vec}^T({\bf I})$ and 
${\bf P}_s({\bf V}_s)={\bf I}-\frac{1}{m}{\rm vec}({\bf V}_s){\rm vec}^T({\bf V}_s^{-1})$, respectively.
Finally, note that for  $s({\bf \Sigma})=|{\bf \Sigma}|^{1/m}$, \eqref{eq:Covariance shape ML} reduces to the simple expression:
\begin{equation}
\label{eq:Covariance shape 3 ML}
{\bf R}_{V_{s,ML}}
=\sigma_1
({\bf I}+{\bf K})({\bf V}_s\otimes {\bf V}_s)
-\frac{2\sigma_1}{m}\ve({\bf V}_s)\ve^T({\bf V}_s).
\end{equation}

The ML estimator of ${\boldsymbol \mu}$ and ${\bf \Sigma}$ are the vector 
$\widehat{\boldsymbol \mu}$ and the symmetric positive definite matrix $\widehat{\bf \Sigma}$ that minimize
the negative log-likelihood function equal from \eqref{eq:pdf RES}
\begin{equation}
\label{eq:log-likelihood function}
L({\boldsymbol \mu},{\bf \Sigma})
=\frac{n}{2}\log|{\bf \Sigma}|
-\sum_{i=1}^n \log \left(g[({\bf x}_i-{\boldsymbol \mu})^T {\bf \Sigma}^{-1}({\bf x}_i-{\boldsymbol \mu})]\right).
\end{equation}
Setting the derivatives of 
$L({\boldsymbol \mu},{\bf \Sigma})$
w.r.t. ${\boldsymbol \mu}$ and ${\bf \Sigma}$ to zero yields the following estimation equations:
\begin{eqnarray}
\label{eq:implicit estimation equations mu}
{\bf 0}
&=&
\frac{1}{n}\sum_{i=1}^n
\varphi[({\bf x}_i-\widehat{\boldsymbol \mu})^T \widehat{\bf \Sigma}^{-1}({\bf x}_i-\widehat{\boldsymbol \mu})]({\bf x}_i-\widehat{\boldsymbol \mu})
\\
\label{eq:implicit estimation equations sigma}
\widehat{\bf \Sigma}
&=&
\frac{1}{n}\sum_{i=1}^n
\varphi[({\bf x}_i-\widehat{\boldsymbol \mu})^T \widehat{\bf \Sigma}^{-1}({\bf x}_i-\widehat{\boldsymbol \mu})]({\bf x}_i-\widehat{\boldsymbol \mu}) ({\bf x}_i-\widehat{\boldsymbol \mu})^T.
\end{eqnarray}
Except for the case where the function $\varphi$ is constant, the set of implicit equations 
\eqref{eq:implicit estimation equations mu}-\eqref{eq:implicit estimation equations sigma} does not guarantee neither the existence nor the uniqueness of the ML estimators $\widehat{\boldsymbol \mu}$ and 
$\widehat{\bf \Sigma}$. Some sufficient conditions to ensure their existence and uniqueness are given in \cite{Kent91}.
Note that for the Gaussian distribution, $\varphi(t)=1$ from \eqref{eq:gauss_cg_dg} and 
\eqref{eq:implicit estimation equations mu}-\eqref{eq:implicit estimation equations sigma}
yield the sample mean $\widehat{\boldsymbol \mu}=\frac{1}{n}\sum_{i=1}^n{\bf x}_i$ and the biased sample covariance matrix
$\widehat{\bf \Sigma}=\frac{1}{n}\sum_{i=1}^n({\bf x}_i-\widehat{\boldsymbol \mu})({\bf x}_i-\widehat{\boldsymbol \mu})^T$, respectively.
%
\subsection{$M$-estimators}
\label{sec:M-estimators}
%
Since the ML estimator may be drastically affected by the presence of outliers
or when the data distribution deviates slightly from the RES distribution of the model, robust estimators have been proposed. Among the different families of robust estimators, in the following we focus our attention on the class of $M$-estimators. As for the ML estimator, an $M$-estimator can be obtained form the minimization of a function on the observation with respect to the parameters of interest. A classical example consists in replacing in the negative log-likelihood \eqref{eq:log-likelihood function},
the function $-\log (g([.])$ by a loss function $\rho(.): \mathbb{R}^+ \mapsto \mathbb{R}$ (generally not related to $g$) to form another objective function. If  $\rho$ is also continuously differentiable with $u(t)\pardef 2 \rho'(t)$, an $M$-estimator is obtained by minimizing this new 
objective function. By replacing $\varphi(t)$ with $u(t)$ and setting again the derivatives of this objective function to zero, 
we obtain similar equations to \eqref{eq:implicit estimation equations mu} and 
\eqref{eq:implicit estimation equations sigma}. 
Some sufficient conditions are also given in \cite{Kent91} to ensure the existence and uniqueness of this $M$-estimator.

The seminal paper by Maronna \cite{maronna1976robust} defined a more general class of $M$-estimator by replacing $\varphi(.)$
by the weight functions $u_1([.]^{1/2})$ and $u_2(.)$ in \eqref{eq:implicit estimation equations mu} and 
\eqref{eq:implicit estimation equations sigma}, respectively, i,e., 
\begin{eqnarray}
\label{eq:Maronna implicit estimation equations mu}
{\bf 0}
&=&
\frac{1}{n}\sum_{i=1}^n
u_1\left([({\bf x}_i-\widehat{\boldsymbol \mu})^T \widehat{\bf \Sigma}^{-1}({\bf x}_i-\widehat{\boldsymbol \mu})]^{1/2}\right)({\bf x}_i-\widehat{\boldsymbol \mu})
\\
\label{eq:Maronna implicit estimation equations sigma}
\widehat{\bf \Sigma}
&=&
\frac{1}{n}\sum_{i=1}^n
u_2\left(({\bf x}_i-\widehat{\boldsymbol \mu})^T \widehat{\bf \Sigma}^{-1}({\bf x}_i-\widehat{\boldsymbol \mu})\right)({\bf x}_i-\widehat{\boldsymbol \mu}) ({\bf x}_i-\widehat{\boldsymbol \mu})^T.
\end{eqnarray}
The functions $u_1([.]^{1/2})$ and $u_2(.): \mathbb{R}^+ \mapsto \mathbb{R}$ need not be the same function and hence the more general 
$M$-estimators need not be related to a minimization problem.
Under sufficient conditions (called Maronna conditions), it is proved in \cite[Th. 4]{maronna1976robust}, 
the existence and uniqueness of solution 
$(\widehat{\boldsymbol \mu},\widehat{\bf \Sigma})$ of 
\eqref{eq:Maronna implicit estimation equations mu},
\eqref{eq:Maronna implicit estimation equations sigma}
but also \cite[Th. 2]{maronna1976robust}  of
$({\bf t},{\bf V})$ solution of 
\begin{eqnarray}
\label{eq:Esperance Maronna implicit equations mu}
{\bf 0}
&=&
\e\left[
u_1\left([({\bf x}_1-{\bf t})^T {\bf V}^{-1}({\bf x}_1-{\bf t})]^{1/2}\right)({\bf x}_1-{\bf t})\right]
\\
\label{eq:Esperance Maronna implicit equations sigma}
{\bf V}
&=&
\e\left[
u_2\left(({\bf x}_1-{\bf t})^T {\bf V}^{-1}({\bf x}_1-{\bf t})\right)({\bf x}_1-{\bf t}) ({\bf x}_1-{\bf t})^T\right].
\end{eqnarray}
Sufficient conditions are also given in \cite[Th. 5]{maronna1976robust} to ensure the strong consistency 
of the estimate $(\widehat{\boldsymbol \mu},\widehat{\bf \Sigma})$ solution of
\eqref{eq:Maronna implicit estimation equations mu},
\eqref{eq:Maronna implicit estimation equations sigma} to the solution 
$({\bf t},{\bf V})$ of 
\eqref{eq:Esperance Maronna implicit equations mu}, 
\eqref{eq:Esperance Maronna implicit equations sigma} with ${\bf t}={\boldsymbol \mu}$
and ${\bf V}=\sigma^{-1}{\bf \Sigma}$, where $\sigma$ is the unique solution of 
$\e[\sigma \mathcal{Q} u_2(\sigma\mathcal{Q})]=m$ \cite[Appendix 3]{tyler1982radial}.

Using a general result on $M$-estimators given in \cite[Sec. 4]{Huber67},
Maronna proved in\cite[Th. 6]{maronna1976robust} the asymptotic gaussianity of 
$(\widehat{\boldsymbol \mu},\widehat{\bf \Sigma})$, where $\widehat{\boldsymbol \mu}$ and $\widehat{\bf \Sigma}$ are asymptotically independent and where only the covariance of the asymptotic distribution of 
$\widehat{\boldsymbol \mu}$ was specified.
Then, using the affine invariance property of any $M$-estimators and the general structure of the covariance of radial random matrices, the
covariance of the asymptotic distribution of 
$\widehat{\bf \Sigma}$ was specified in \cite[Appendix 2]{tyler1982radial} to get:
\begin{equation}
\label{eq:th M asymptotic}
\sqrt{n}
\left(\begin{array}{c}
\widehat{\boldsymbol \mu}-{\boldsymbol \mu} \\
\ve(\widehat{\bf \Sigma})-\ve({\bf \sigma^{-1}\Sigma})
\end{array}
\right)
 \rightarrow_d \mathbb{R}N_{m+m^2}
\left(
\left(\begin{array}{c}
{\bf 0}\\
{\bf 0} 
\end{array}\right),
\left(
\begin{array}{cc}
{\bf R}_{{\mu}_{M}}& {\bf 0}\\
{\bf 0} & {\bf R}_{{\Sigma}_{M}}\\
\end{array}\right)
\right),
\end{equation}
where
\begin{equation}
\label{eq:Rmu M}
{\bf R}_{{\mu}_{M}}
=\frac{\alpha}{\beta^2}{\bf V}
=\frac{\alpha \sigma^{-1}}{\beta^2}{\bf \Sigma}
\end{equation}
where
$\alpha=\frac{1}{m}\e[\psi_1^2(\sqrt{{\sigma}\mathcal{Q}})]$
and
$\beta=\e[(1-m^{-1})u_1(\sqrt{{\sigma}\mathcal{Q}})+m^{-1}\psi'_1(\sqrt{{\sigma}\mathcal{Q}})]$
with $\psi_1(t)\pardef tu_1(t)$, and
\begin{equation}
\label{eq:Rsigmauu M}
{\bf R}_{{\Sigma}_{M}}
=\sigma_1({\bf I}+{\bf K})({\bf V}\otimes{\bf V})
+\sigma_2\ve({\bf V})\ve^T({\bf V}),
\end{equation}
where
$\sigma_1
=\frac{(m+2)^2a_1}{(2a_2+m)^2}$
and
$\sigma_2
=a_2^{-1}
\left[
(a_1-1)-2(a_2-1)a_1\frac{m+(m+4)a_2}{(2a_2+m)^2}
\right]$
with
$a_1=\frac{\e[\psi_2^2({\sigma}\mathcal{Q})]}
{m(m+2)}$
and
$a_2=\frac{\e[{\sigma}\mathcal{Q}{\psi_2}^{'}({\sigma}\mathcal{Q})]}{m}$
where $\psi_2(t) \pardef t u_2(t)$ and $\psi_2^{'}(t) \pardef \frac{d\psi_2(t)}{dt}$.
Furthermore, using the general structure of the covariance of radial random matrices, 
it is proved in \cite[Th. 1]{tyler1982radial} that
\begin{equation}
\label{eq:sigma 1 et 2}
\sigma_2 \ge
-\frac{2}{m}\sigma_1.
\end{equation}
Finally note  that for $u_1(t)=u_2(t^2)=\varphi(t^2)$, the $M$-estimate reduces to the ML estimate for which $\sigma=1$ and \eqref{eq:th M asymptotic} reduce to \eqref{eq:th ML asymptotic}.

To the best of our knowledge, when both ${\boldsymbol \mu}$ and ${\bf \Sigma}$ are unknown parameters, 
it does not exist in the literature an analysis of sufficient conditions under which an algorithm would jointly converge toward the solution of the implicit equations 
\eqref{eq:Maronna implicit estimation equations mu}, \eqref{eq:Maronna implicit estimation equations sigma}.
However, when ${\boldsymbol \mu}$ is known, an algorithm which is essentially a fixed-point algorithm with a scale adjustment made at each
iteration, with less severe conditions than that required in \cite{maronna1976robust} was presented in \cite{tyler1988some}.
Without loss of generality, ${\boldsymbol \mu}$ can be taken to ${\bf 0}$ and the following algorithm is proved in \cite{tyler1988some} to converge to the unique solution of 
$\widehat{\bf \Sigma}
=
\frac{1}{n}\sum_{i=1}^n
u_2[{\bf x}_i^T \widehat{\bf \Sigma}^{-1}{\bf x}_i]{\bf x}_i{\bf x}_i^T$:
\begin{equation}
\label{eq:algo M estimate sigma}
{\bf \Sigma}_{0}=\frac{1}{n}\sum_{i=1}^n{\bf x}_i{\bf x}_i^T, \ \ 
{\bf \Sigma}_{k+1}
=
\frac{1}{n}\sum_{i=1}^n
u_2[c_k{\bf x}_i^T {\bf \Sigma}_k^{-1}{\bf x}_i]{\bf x}_i{\bf x}_i^T,
\end{equation}
$c_k$ being the unique positive scalar satisfying
$\frac{1}{n}\sum_{i=1}^n
\psi_2[c_k{\bf x}_i^T {\bf \Sigma}_k^{-1}{\bf x}_i]=m$ with $\psi_2(t)\pardef t u_2(t)$.
%
\subsection{Tyler's $M$-estimator}
\label{sec:Tyler's-estimator}
%
When ${\boldsymbol \mu}$ is known (equal to ${\bf 0}$ without loss of generality), the $M$-estimator proposed by 
Tyler \cite{tyler1987distribution} has become a very popular robust scatter estimator in the signal processing literature. This $M$-estimator is defined by its weight function $u_2(t)=\frac{m}{t}$ associated with the loss function $\rho(t)=\frac{m}{2}\log t$, leading to the following objective function
\begin{equation}
\label{eq:Tyler objective function}
L_T({\bf \Sigma})
=\frac{n}{2}\log|{\bf \Sigma}|
+\frac{m}{2}\sum_{i=1}^n \log ({\bf x}_i^T {\bf \Sigma}^{-1}{\bf x}_i)
\end{equation}
that is minimized by $\widehat{\bf \Sigma}$, solution of the implicit equation
\begin{equation}
\label{eq: implicit equation Tyler}
\widehat{\bf \Sigma}
=
\frac{m}{n}\sum_{i=1}^n
\frac{{\bf x}_i{\bf x}_i^T}
{{\bf x}_i^T\widehat{\bf \Sigma}^{-1}{\bf x}_i}.
\end{equation}
Note that $L_T(c^2{\bf \Sigma})= L_T({\bf \Sigma})+a$ where $a$ does not depend on ${\bf \Sigma}$. Consequently, if
$\widehat{\bf \Sigma}$ is solution of \eqref{eq: implicit equation Tyler}, $c^2\widehat{\bf \Sigma}$ is it too.
The existence proof of solution of \eqref{eq: implicit equation Tyler} given in Maronna \cite{maronna1976robust} does not apply here, but
existence and uniqueness (up to a multiplicative constant) under continuous RES distribution in  
$\mathbb{R}^m$ is proved in 
\cite{tyler1987distribution}, by showing that it is the limiting point of the following specific fixed point algorithm: with ${\bf \Sigma}_0$ arbitrary symmetric positive definite matrix,
\begin{equation}
\label{eq: algorithme Tyler}
{\bf \Sigma}'_{k+1}
=
\frac{m}{n}\sum_{i=1}^n
\frac{{\bf x}_i{\bf x}_i^T}
{{\bf x}_i^T{\bf \Sigma}_k^{-1}{\bf x}_i},
\ \ 
{\bf \Sigma}_{k+1}
=\frac{m}{\tra({\bf \Sigma}'_{k+1})}{\bf \Sigma}'_{k+1}.
\end{equation}
Under this same condition, it is also proved in \cite{tyler1987distribution} that the solution of 
\eqref{eq: implicit equation Tyler} standardized so that $\tra(\widehat{\bf \Sigma})=m$ is strongly consistent to the symmetric positive definite matrix ${\bf \Sigma}$ which is also solution of 
\begin{equation}
\label{eq: limite Tyler estimate}
{\bf \Sigma}
=
m\e\left[
\frac{{\bf x}_1{\bf x}_1^T}
{{\bf x}_1^T{\bf \Sigma}^{-1}{\bf x}_1}\right]
\ \ \mbox{with}
\ \ \tra({\bf \Sigma})=m.
\end{equation}
Furthermore, under continuous RES distribution in $\mathbb{R}^m$, it is proved in \cite{tyler1987distribution} that 
the solution $\widehat{\bf \Sigma}$ of
\eqref{eq: implicit equation Tyler} constrained to $\tra(\widehat{\bf \Sigma})=m$  is asymptotically Gaussian distributed, i.e.,
\begin{equation}
\label{eq: Tyler asymptotique}
\sqrt{n}\left(
\ve(\widehat{\bf \Sigma})-\ve({\bf \Sigma})
\right) \rightarrow_d \mathbb{R}N_{m^2}\left({\bf 0},{\bf R}_{{\Sigma}_{Ty}}  \right)
\end{equation}
with
\begin{equation}
\label{eq:asymptotic covariance Tyler}
{\bf R}_{{\Sigma}_{Ty}}
=\left(1+\frac{2}{m}\right)({\bf I}+{\bf K})({\bf \Sigma}\otimes{\bf \Sigma})
-\frac{2}{m}\left(1+\frac{2}{m}\right) \ve({\bf \Sigma})\ve^T({\bf \Sigma}).
\end{equation}
Comparing \eqref{eq:Rsigmauu M} to \eqref{eq:asymptotic covariance Tyler}, we see that
$\sigma_2=-\frac{2}{m}\sigma_1$ and therefore we have equality in
\eqref{eq:sigma 1 et 2}.

Note that when ${\bf x}_i$ are SIRV distributed, i.e., 
${\bf x}_i \sim$ RCG$_m({\bf 0}, {\bf \Sigma}, F_{{\tau}_i})$ where here
$\tau_i$ are assumed deterministic and unknown, the ML estimate of ${\bf \Sigma}$
coincides with Tyler's $M$-estimate \eqref{eq: implicit equation Tyler}. 
This estimate was introduced in \cite{Gini1997} and existence and uniqueness were studied in 
\cite{PCOFL08}.
Tyler's $M$-estimator enjoys four interesting properties:
\begin{itemize}
\item[$\bullet$]
Objective function \eqref{eq:Tyler objective function} is also the negative log-likelihood of i.i.d. data from an ACG$_m({\bf 0},{\bf \Sigma})$ distribution (see \eqref{eq:ACG_pdf}). Consequently Tyler's $M$-estimator of 
${\bf \Sigma}$
is the ML estimator $\widehat{\bf \Sigma}$ (satisfying the constraint $\tra(\widehat{\bf \Sigma})=m$) under the ACG$_m({\bf 0},{\bf \Sigma})$ distribution.
\item[$\bullet$]
Replacing ${\bf x}_i$ by $\frac{{\bf x}_i}{\|{\bf x}_i\|}$ in \eqref{eq: implicit equation Tyler} does not affect the solution of \eqref{eq: implicit equation Tyler}. Since $\frac{{\bf x}_i}{\|{\bf x}_i\|}$ is 
ACG$_m({\bf 0},{\bf \Sigma})$ distributed for arbitrary RES$_m({\bf 0},{\bf \Sigma},g)$ distribution, the distribution of Tyler's $M$ estimator $\widehat{\bf \Sigma}$ 
does not depend on the density generator of this RES distribution.
This is the reason why Tyler referred to his estimate as a distribution-free estimator \cite{tyler1987distribution}.
\item[$\bullet$]
Suppose ${\bf x}_1,..,{\bf x}_i,..,{\bf x}_n$ are independent where ${\bf x}_i \sim {\rm RES}_m({\bf 0},c_i^2{\bf \Sigma},g_i)$ with unknown parameters $c_i^2$ and arbitrary unknown density generators $g_i$, then the ML of ${\bf \Sigma}$ corresponds to Tyler's $M$-estimate
\cite[Th. 1]{ollila2012distribution}.
\item[$\bullet$]
Tyler's $M$-estimator can be considered as the most robust estimator of the scatter matrix of a RES distribution in the sense of minimizing the maximum asymptotic covariance w.r.t. the generator density (see \cite[Remark 3.1]{tyler1987distribution}, and \cite[Th. 1]{tyler1983robustness}.
\end{itemize}

Also, note that Tyler's $M$-estimate has been extensively studied in the statistics literature
(e.g., \cite{KING1994} has specified the existence proof of solution of \eqref{eq: implicit equation Tyler} and \cite{FM2020} have proved that Tyler's iterative procedure 
\eqref{eq: algorithme Tyler} has a linear convergence rate).

We finally note that many complementary results have been carried out on RES distributions
(see e.g., \cite{PBTB2013} for ML estimate of GG scatter, \cite{SW2014} for Tyler's estimate of structured scatter). 
Furthermore, most of the definitions and results presented in this section have been extended to CES distributions with many new results (see, e.g., 
\cite{PCOFL08,PFOL2008,CP2008,MPFO2013,DP2018,MREBF18,MREKBF2019}).
%
\subsection{Slepian-Bangs formula}
\label{sec:Slepian-Bangs formula}
We use in this Section the real to complex representation to unify the Slepian-Bangs (SB) for RES, C-CES and NC-CES distributed data.
The SB formula provides a convenient way to compute the FIM and thus the Cramer-Rao bound (CRB) on the real-valued parameter ${\boldsymbol \alpha}$ parameterizing and characterizing the couple $({\boldsymbol \mu},{\bf \Sigma})$ of elliptically symmetric distributions.
We omit this dependence in ${\boldsymbol \alpha}$ to simplify the notations.
To derive this formula, it is necessary that the scatter matrix ${\bf \Sigma}$ is not singular and the second-order moments of ${\bf x}$ are finite,
which is equivalent to the first-order moments of $\mathcal{Q}$ being finite.
Furthermore, to avoid the ambiguity between the scatter matrix and the density generator,
we either assume that ${\bf \Sigma}=\cov({\bf x})$ 
or that there is a scale constraint on ${\bf \Sigma}$.
This formula has been derived for
the real Gaussian distribution in \cite{S54} and \cite{B71}, then extended to the circular complex Gaussian and non-circular Gaussian case in \cite{SM05} and \cite{DA04}, respectively.
This formula has been extended to C-CES distributions in \cite{BA13} and \cite{GG13a}, and recently to NC-CES distributions \cite{AD19}.
In all these scenarios, the density generator is assumed to be perfectly known.
For RES distributed ${\bf x}$, all the steps of the proof of the SB formula for C-CES distributions given in \cite{BA13} 
apply, and we get the following structured matrix SB formula:
\begin{eqnarray}
\nonumber
{\rm CRB}^{-1}({\boldsymbol \alpha})
&=&
a_0\frac{d{\boldsymbol \mu}^T}{d {\boldsymbol \alpha}^T}{\bf \Sigma}^{-1}\frac{d{\boldsymbol \mu}}{d {\boldsymbol \alpha}^T}
\\
\nonumber
&+&
\left(
\frac{d \ve({\bf \Sigma})}{d{\boldsymbol \alpha}^T}\right)^T
\left(a_1 ({{\bf \Sigma}}^{-T}\otimes{\bf \Sigma}^{-1})+
a_2 \ve({\bf \Sigma}^{-1})\ve^T({\bf \Sigma}^{-1})
\right)
\\
\label{eq:structure FIM}
&&
\hspace{6cm}
\left(\frac{d \ve({\bf \Sigma})}{d {\boldsymbol \alpha}^T}
\right),
\end{eqnarray}
where $a_0=\xi_{r,1,m}$, $a_1=\frac{1}{2}\xi_{r,2,m}$ and $a_2=\frac{1}{4}(\xi_{r,2,m}-1)$ with
\begin{equation}
\label{eq:def xi1r and x2r}
\xi_{r,1,m}\pardef\frac{\e[Q\varphi_r^2(Q)]}{m}
\ \ \mbox{and}\ \
\xi_{r,2,m}\pardef\frac{\e[Q^2\varphi_r^2(Q)]}{m(m+2)},
\end{equation}
where $Q \pardef \mathcal{Q}_{r,m}$ and
$\varphi_r(t)\pardef  -\frac{2}{g_{r,m}(t)}\frac{d g_{r,m}(t)}{dt}$.

Note that because 
$(c^2{\bf \Sigma}, c^{-2}Q,c^{m}g_{r,m}(.\ c^2), c^2\varphi_r(.\ c^2))$ gives the same RES distribution,
the coefficient $\xi_{r,2,m}$ is free of scale ambiguity in contrast to $\xi_{r,1,m}$ which depends on the scale factor. This is consistent with eq. \eqref{eq:structure FIM}.

This SB formula allows us to directly deduce those of  NC-CES distributed data
obtained thanks to the real to complex representation. This SB formula is similarly structured where
$\boldsymbol \mu$, $\bf \Sigma$, 
$\frac{d{\boldsymbol \mu}^T}{d {\boldsymbol \alpha}^T}$, $\left(\frac{d \ve({\bf \Sigma})}{d{\boldsymbol \alpha}^T}\right)^T$
and $\ve^T({\bf \Sigma}^{-1})$
in \eqref{eq:structure FIM} are replaced by 
$\widetilde{\boldsymbol \mu}$, $\widetilde{\bf \Sigma}$, 
$\frac{d{\widetilde{\boldsymbol \mu}}^H}{d {\boldsymbol \alpha}^T}$, $\left(\frac{d \ve(\widetilde{\bf \Sigma})}{d{\boldsymbol \alpha}^T}\right)^H$ and
 $\ve^H(\widetilde{\bf \Sigma}^{-1})$, respectively, where $a_0=\xi_{c,1,m}$, $a_1=\frac{\xi_{c,2,m}}{2}$ and 
$a_2=\frac{1}{4}(\xi_{c,2,m}-1)$
with
\begin{equation}
\label{eq:def xi1c and x2c}
\xi_{c,1,m}\pardef \frac{\e[Q\varphi_c^2(Q)]}{m}
\ \ \mbox{and}\ \
\xi_{c,2,m}\pardef \frac{\e[Q^2\varphi_c^2(Q)]}{m(m+1)},
\end{equation}
where $Q \pardef \mathcal{Q}_{c,m}$ and $\varphi_c(t)\pardef  -\frac{1}{g_{c,m}(t)}\frac{d g_{c,m}(t)}{dt}$.
On the other hand, the SB formulas for C-CES distributed data can be deduced directly  by replacing  $\widetilde{\bf \Sigma}$ by 
$\footnotesize{\left(\begin{array}{cc}
{\bf \Sigma}& {\bf 0}\\
{\bf 0} &  {\bf \Sigma}^* \\
\end{array}
\right)}$, yielding the SB formulas proved in
\cite{BA13} and \cite{GG13a}, which are also similarly structured with $a_0=2\xi_{c,1,m}$, $a_1=\xi_{c,2,m}$ and $a_2=\xi_{c,2,m}-1$, 
where
$\frac{d{\boldsymbol \mu}^T}{d {\boldsymbol \alpha}^T}{\bf \Sigma}^{-1}\frac{d{\boldsymbol \mu}}{d {\boldsymbol \alpha}^T}$,
$\left(\frac{d \ve({\bf \Sigma})}{d{\boldsymbol \alpha}^T}\right)^T$
and $\ve^T({\bf \Sigma}^{-1})$ are replaced in 
by
${\rm Re}\left(\frac{d{\boldsymbol \mu}^H}{d {\boldsymbol \alpha}^T}{\bf \Sigma}^{-1}\frac{d{\boldsymbol \mu}}{d {\boldsymbol \alpha}^T}\right)$,
$\left(\frac{d \ve({\bf \Sigma})}{d{\boldsymbol \alpha}^T}\right)^H$ and $\ve^H({\bf \Sigma}^{-1})$,
respectively.

Note that $(a_0,a_1,a_2)$ reduces to $(1,1/2,0)$, $(1,1/2,0)$ and $(2,1,0)$ for real, complex noncircular and complex circular Gaussian distributions, respectively.
Note also that the decoupling between the parameters $\boldsymbol{\alpha}_1$ of $\boldsymbol{\mu}$ and the parameters $\boldsymbol{\alpha}_2$ of 
$\boldsymbol{\Sigma}$ when $\boldsymbol{\mu}$ and $\boldsymbol{\Sigma}$ have no parameters in common for Gaussian distributions, extends to any 
elliptically symmetric distributions
with  
\begin{equation}
\label{eq:CRB alpha1}
{\rm CRB}({\boldsymbol \alpha}_1)
=
\left(a_0\frac{d{\boldsymbol \mu}^T}{d {\boldsymbol \alpha_1}^T}{\bf \Sigma}^{-1}\frac{d{\boldsymbol \mu}}{d {\boldsymbol \alpha_1}^T}\right)^{-1} 
\end{equation}
and
\begin{eqnarray}
\nonumber
{\rm CRB}({\boldsymbol \alpha}_2)
&=&
\left(\left(
\frac{d \ve({\bf \Sigma})}{d{\boldsymbol \alpha_2}^T}\right)^T\
\left(a_1 ({{\bf \Sigma}}^{-T}\otimes{\bf \Sigma}^{-1})+
a_2 \ve({\bf \Sigma}^{-1})\ve^T({\bf \Sigma}^{-1})
\right)\right.
\\
\label{eq:CRB alpha2}
&&
\hspace{6cm}
\left.
\left(\frac{d \ve({\bf \Sigma})}{d {\boldsymbol \alpha_2}^T}
\right)\right)^{-1}.
\end{eqnarray}

The coefficients $\xi_{c,1,m}$ and $\xi_{c,2,m}$ are deduced from $\xi_{r,1,m}$ and $\xi_{r,2,m}$ by the general relations
\begin{equation}
\label{exi r etc}
\xi_{c,1,m}=\xi_{r,1,2m} 
\ \ \mbox{and}\ \ \xi_{c,2,m}=\xi_{r,2,2m}.
\end{equation}
For example, the coefficients $\xi_{c,1,m}$ and $\xi_{c,2,m}$ have been calculated for complex Student's $t$ and complex generalized Gaussian distributions in \cite{BA13} and \cite{GG13a}. They are given respectively by:
\begin{eqnarray}
\label{eq:xi1 et xi2 student}
\xi_{c,1,m}
&=&
\frac{\nu/2}{((\nu/2)-1)}\frac{(\nu/2)+m}{((\nu/2)+m+1)}
\ \mbox{and} \ \
\xi_{c,2,m}=\frac{(\nu/2)+m}{(\nu/2)+m+1},
\\
\label{eq:xi1 et xi2 CGG}
\xi_{c,1,m}
&=&
\frac{\Gamma(2+\frac{m-1}{s}) \Gamma(\frac{m+1}{s})}{(\Gamma(1+\frac{m}{s}))^2}
 \hspace{1.2cm} \mbox{and} \ \
\xi_{c,2,m}
=\frac{m+s}{m+1}.
\end{eqnarray}
Finally note that the SB formula for elliptically symmetric distributions has been extended
when the data model is misspecified by the parametric probabilistic model 
\cite{For_SCRB}, and when the density generator is considered as an infinite-dimensional nuisance parameter \cite{MFEYZB18} or
parameterized by a nuisance parameter \cite{AD2023}.
%
\section{Conclusion}
\label{sec:conclusion} 
%
The aim of this chapter was to provide a short overview of the main properties of elliptically symmetric distributions and it can be used as a background for all the other chapters in this book. There is no claim to completeness in the material presented here. The reader interested in deeper discussions and investigations on specific aspects may find the references list useful to this goal. As a last concluding remark, we would like to highlight our choice to focus mainly on the RES distributions. As explained in the chapter, this class can be considered as the most general one since it encompasses the C-CES and NC-CES distributions as special cases. Some effort has been then put into showing explicitly the mapping between the RES class and all its sub-class. We hope that this chapter may represent a reference for the reader helping him to not get lost while going down his path through this book.
%
\bibliographystyle{IEEEtran}
\bibliography{bibliojp}
\end{document}